\documentclass[a4paper,10pt]{amsart}
\usepackage{amssymb}
\usepackage{amsmath}
\usepackage{epsfig}
\usepackage[all]{xy}
\usepackage{comment}

\newtheorem{thm}{Theorem}
\newtheorem{cor}[thm]{Corollary}
\newtheorem{lem}[thm]{Lemma}
\newtheorem{prop}[thm]{Proposition}
\newtheorem{defn}[thm]{Definition}
\newtheorem{conjecture}[thm]{Conjecture}

\newtheorem{theorem-question}[thm]{Theorem-Question}

\newenvironment{exafont}{\begin{bf}}{\end{bf}}

\newenvironment{remark}{\vspace{0.3cm}\par\noindent\refstepcounter{thm}\begin{exafont}Remark \thethm\end{exafont}\hspace{\labelsep}}{\vspace{0.3cm}\par}

\newcommand{\cS}{\mathcal{S}}
\newcommand{\cC}{\mathcal{C}}
\newcommand{\cN}{\mathcal{N}}
\newcommand{\cG}{\mathcal{G}}
\newcommand{\cJ}{\mathcal{J}}
\newcommand{\cL}{\mathcal{L}}
\newcommand{\cQ}{\mathcal{Q}}
\newcommand{\cK}{\mathcal{K}}

\newcommand{\la}{\lambda}
\newcommand{\ZZ}{\mathbb{Z}}
\newcommand{\eps}{\varepsilon}
\DeclareMathOperator{\End}{End}
\DeclareMathOperator{\Hom}{Hom}
\DeclareMathOperator{\Ext}{Ext}
\DeclareMathOperator{\add}{add}
\DeclareMathOperator{\soc}{soc}
\DeclareMathOperator{\im}{im}
\DeclareMathOperator{\F}{\mathcal{F}}

\DeclareMathOperator{\ml}{-mod}
\DeclareMathOperator{\mr}{mod-}
\DeclareMathOperator{\cat}{-cat}
\DeclareMathOperator{\sml}{-\underline{mod}}
\DeclareMathOperator{\projl}{-proj}

\begin{document}

\baselineskip = 15pt

\title{\large Rational representations of $GL_2$.} 

\author{Vanessa Miemietz and Will Turner}

\thanks{The first author acknowledges support from Leverhulme.}

\maketitle

\begin{abstract}
Let $F$ be an algebraically closed field of characteristic $p$.
We fashion an infinite dimensional basic algebra $\underleftarrow{\cC}_p(F)$, 
with a transparent combinatorial structure, 
which we expect to control the rational representation theory of $GL_2(F)$.
\end{abstract}

\section{Introduction}

In any first course on representation theory,
the students will become familiar with 
representations of the algebraic group $GL_2(\mathbb{C})$, or with those of some close relation of this group. 
It is perhaps surprising therefore, over a century after the birth of group representation theory,
 that anything remains to be said about $GL_2$.
However, development in the modular theory has been much slower than in
characteristic zero and even of the smallest cases no full understanding
has yet been reached.
In this article, we wish to pursuade the reader that there is structure underlying the rational representation theory
of $GL_2$ over a field of positive characteristic, 
as simple as the structure appearing in characteristic zero, although quite different in nature.

Of course, even in positive characteristic,
the usual hare-headed questions about $GL_2$-modules were answered long ago:
irreducibles are parametrised by elements of the dominant region of the weight lattice,
and have realisations as tensor products of 
Frobenius twists of socles of 
symmetric powers of the natural representation in small degrees, 
and powers of the determinant representation.  
However, the situation is more delicate than these easy truths imply. 
There are homological interactions between irreducible modules, 
and for a deeper understanding one ought to contemplate the manner in which these interactions occur.
This is the concern of our paper.

\bigskip

We shall be more precise.
Let $A$ be an algebra with a self-dual bimodule $T$.
Let $B$ be the algebra whose category of ungraded representations
is equivalent to the category of graded representations of the trivial extension algebra of $A$ by $T$.
Let $C$ be the trivial extension algebra of $B$ by its dual.
Modulo the infinite dimensionality of $C$, we have a map
$${\cC} \circlearrowright \{ \textrm{algebras with a self-dual bimodule} \},$$
which takes an algebra $A$, with an $A$-$A$ bimodule $T$, such that 
$_AT_A \cong  {_AT^*_A}$, to a symmetric algebra $C$.
The self-dual bimodule corresponding to $C$ is the
regular bimodule $_CC_C$.

For every $n \in \mathbb{N}$, there is a localisation
$${\cC}_n \circlearrowright \{ \textrm{algebras with a self-dual bimodule} \}$$
of ${\cC}$.
There is a canonical epimorphism
$$A \twoheadleftarrow {\cC}_n(A).$$
Taking the inverse limit of the sequence
$$A \twoheadleftarrow {\cC}_n(A)  \twoheadleftarrow {\cC}_n({\cC}_n(A))
\twoheadleftarrow {\cC}_n({\cC}_n({\cC}_n(A)))  \twoheadleftarrow...,$$
we obtain an algebra $\underleftarrow{\cC}_n(A)$.

\bigskip

Let ${\cS}(2) = \bigoplus_{r \geq 0} {\cS}(2,r)$ be the Schur algebra associated to 
$GL_2$, defined over an algebraically closed field $F$ of characteristic $p>0$ \cite{Green}.
There is a sequence of natural surjections
$${\cS}(2,r) \twoheadleftarrow {\cS}(2,r+2) 
\twoheadleftarrow {\cS}(2,r+4) \twoheadleftarrow {\cS}(2,r+6) \twoheadleftarrow...$$
Let ${\cS}(2,\underline{r})$ 
be the inverse limit of this directed sequence of algebra 
epimorphisms. 
The category of rational representations of $GL_2(F)$ is equivalent to the category of 
finite dimensional representations
of the direct sum $\bigoplus_{r \in \mathbb{Z}} {\cS}(2, \underline{r})$.

In the sequel, we define a certain filtration on ${\cS}(2,r)$, refining the radical filtration, and denote by
${\cG}(2,r)$ the graded ring associated to this filtration.
There is a compatible sequence of surjections
$${\cG}(2,r) \twoheadleftarrow {\cG}(2,r+2) 
\twoheadleftarrow {\cG}(2,r+4) \twoheadleftarrow {\cG}(2,r+6) \twoheadleftarrow...$$
Let ${\cG}(2,\underline{r})$ 
be the inverse limit of this directed sequence of algebra epimorphisms. 
Our main result is the following:

\begin{thm} \label{Morita}
Every block of ${\cG}(2, \underline{r})$ 
is Morita equivalent to $\underleftarrow{\cC}_p(F)$.
\end{thm}

Our proof of Theorem \ref{Morita} is inductive. 
We apply results of K. Erdmann, A. Henke, and S. Koenig concerning ${\cS}(2,r)$ (\cite{EH:02}, \cite{HK:02}), 
to prove that certain Ringel self-dual blocks of ${\cG}(2, r)$ are equivalent to
${\cC}_p^d(F)$, for some $d$.
Since every block of ${\cG}(2, r)$ is a quotient of such a Ringel self-dual block,
the theorem follows.

In fact, we prove a rather stronger statement.
Let $\cS$ be a Ringel self-dual Schur algebra ${\cS}(2,r)$. 
We demonstrate the existence of a filtration by ideals,
$$\cS \supset {\cN} \supset {\cN}^2 \supset 0,$$ 
whose associated graded ring is Morita equivalent to ${\cC}_a(A) \oplus F^{\oplus m}$, 
where $A$ is a smaller Ringel self-dual Schur algebra ${\cS}(2,s)$, where $2 \leq a \leq p$,
and where $m$ is some multiplicity. 

In an earlier chapter, we give careful definitions of  
$B,C$, and ${\cC}_p$, 
and prove that under favourable conditions, they respect the quasi-heredity condition.
 
\bigskip

This is all very pleasing, but we believe more to be true. 
We predict that in fact ${\cS}(2,r) \cong {\cG}(2,r)$, for all $r$, and therefore
${\cS}(2, \underline{r}) \cong {\cG}(2, \underline{r})$.
In other words, we expect the following to be true:

\begin{conjecture} \label{wishful}
Every block of rational representations of $GL_2(F)$ is equivalent to 
$\underleftarrow{\cC}_p(F) \ml$.
\end{conjecture}

In the final chapter of the paper, we consider this possibility in more detail. 
We demonstrate that the main obstacle to a proof by induction
is a familiar one in modular representation theory: the lifting of a stable equivalence.

In his inductive approach to M. Brou\'e's abelian defect group conjecture \cite{Broue},
R. Rouquier has established that the main difficulty 
is the lifting of a stable equivalence to a derived equivalence \cite{Rouquier}.
In our microcosm, we give a similar inductive strategy 
to prove that ${\cS}(2,r) \cong {\cG}(2,r)$.
We define a pair of infinite dimensional self-injective algebras, ${\cL}_1$ and ${\cL}_2$,
and prove the existence of a stable equivalence between these, 
sending simple modules to simple modules.  
If one could lift this stable equivalence to a Morita equivalence, 
an isomorphism ${\cS}(2,\underline{r}) \cong {\cG}(2,\underline{r})$ would follow.

\bigskip

There are ramifications for the Temperley-Lieb algebra, which we briefly mention here.
If $F$ has characteristic $p>2$, 
then an $r$-fold tensor product of the natural two dimensional $GL_2(F)$-module
is a full tilting module for ${\cS}(2,r)$.
Therefore its endomorphism ring, 
known as the Temperley-Lieb algebra $TL_r$, is the Ringel dual of ${\cS}(2,r)$.
We have directed sequences of embeddings of idempotent subalgebras,
$$TL_r \hookrightarrow TL_{r+2} \hookrightarrow TL_{r+4} \hookrightarrow TL_{r+6} \hookrightarrow...$$
$$A \hookrightarrow {\cC}_n(A) \hookrightarrow {\cC}_n({\cC}_n(A))
\hookrightarrow {\cC}_n({\cC}_n({\cC}_n(A))) \hookrightarrow ... .$$
Let $TL_{\underline{r}}, \underrightarrow{\cC}_n(A)$ denote the direct limits of these
sequences of algebra monomorphisms. We expect any block of
$TL_{\underline{r}}$ to be Morita equivalent to $\underrightarrow{\cC}_p(F)$.

\bigskip

\emph{Acknowledgements:}
We have benefited from conversations with Anne Henke and Karin Erdmann concerning Schur algebras, 
with Joe Chuang, and with Rapha\"el Rouquier.
We thank the EPSRC for financial support.

\section{Setup}

Throughout this paper, $F$ will be a field and $A$ an $F$-algebra. 
We denote by $A \ml$ the category of finite dimensional left $A$-modules, and by $\mr A$ 
the category of finite dimensional right $A$-modules.
Given a finite dimensional left/right module $M$, we write the dual of $M$ as $M^* = \Hom_F(M,F)$, 
a right/left module.

We write $A \projl$ for the category of finite dimensional projective left 
$A$-modules.
Given a collection $X \subset A \ml$, we denote by ${\F}(X)$ the category
of modules filtered by objects in $X$.
Let ${\cJ}(A)$ denote the Jacobson radical of $A$.
 
We suppose that $A$ is a locally finite dimensional algebra.
In other words, there exists a set $\Lambda$, 
indexing a set of orthogonal idempotents $\{ e_\lambda \}_{\lambda \in \Lambda}$, 
such that $A \cong \bigoplus_{\lambda,\mu \in \Lambda} e_\lambda A e_\mu$, 
and $e_\lambda A e_\mu$ is finite dimensional.
We assume that $A/{\cJ}(A) = \bigoplus_{\lambda \in \Lambda} M_\lambda$ 
is a direct sum of matrix rings $M_\lambda$ over $F$, where $e_\lambda$ is the unit of $M_\lambda$.
Thus, $\Lambda$ is an indexing set for isomorphism classes of simple $A$-modules. 
By the idempotent decomposition, simple modules have projective covers and injective hulls, 
providing $1-1$-correspondences between isomorphism classes of simples, projectives and injectives.

Now let $\Lambda$ be a poset which is interval-finite 
(i.e. for every $\mu \leq \la \in \Lambda$ the set $\{ \nu | \mu \leq \nu \leq \la\}$ is finite).
 
Recall that $\mr A$
is a highest weight category in the sense of Cline, Parshall and Scott \cite{CPS:88} if,
for every $\lambda \in \Lambda$ there exists an irreducible right module $L^r(\la)$, 
a costandard right module
$\nabla^r(\la)$, which embeds into the injective hull $I^r(\la)$ of $L^r(\la)$,
such that the cokernel of this inclusion is filtered by
$\nabla^r(\mu)$ for $\mu \geq \la$, and $\nabla^r(\la)/ \soc
\nabla^r(\la) $ consists of composition factors $L^r(\nu)$ for
$\nu < \la$.  
Dualizing with respect to $F$, we find this is equivalent to the corresponding
projective indecomposable left modules $P(\la) \in A \ml$ 
having standard filtrations.  
So, for
every $\la \in \Lambda$ there exists a standard module $\Delta(\la)$
and an epimorphism $P(\la) \twoheadrightarrow \Delta(\la)$, the kernel
of which is filtered by modules $\Delta(\mu)$ for $\mu > \la$, and the
kernel of the map $\Delta(\la)\twoheadrightarrow L^l(\la)$ consists of
composition factors of the form $L(\nu)$ for $\nu < \la$.

Let $J \subset \Lambda$ be a nonempty finitely generated ideal.
The subcategory $\mr A [J]$ of objects which only have composition factors
$L(\nu)$ for $\nu \in J$ is a highest weight category, whenever $\mr A$ is a highest weight category
(\cite{CPS:88}, Theorem 3.5). 
Let $A^J  = A/\sum_{\lambda \notin J} Ae_\lambda A$.
Then $A^J$ is a locally finite-dimensional algebra and $\mr A [J] \cong \mr A^J$.

Let $I \subset \Lambda$ be a  nonempty finitely generated coideal  and 
define 
$A_I := \underset{\la, \mu \in I}{\bigoplus} e_\la A e_\mu.$

\begin{lem}
If $\mr A$ is a highest weight category, then $\mr A_I$ is a highest weight category.
\end{lem}

\emph{Proof.} 
We construct $\Delta$-filtrations of projectives in $A_I \ml$. Projectives in $A_I \ml$ are of the form 
$P_{A_I}(\lambda) := \Hom_A(\underset{ \mu \in I}{\bigoplus}Ae_\mu ,P_{A}(\la))$. 
We define
$\Delta_{A_{I}}(\lambda):= \Hom_A(\underset{ \mu \in I}{\bigoplus} A_\mu,\Delta_{A}(\la))$.
Since $\Hom_A(\underset{ \mu \in I}{\bigoplus} A_\mu,-)$ is exact we obtain a filtration of $P_{A_I}(\la)$ respecting the necessary conditions on orders.
$\square$

\bigskip

Let us define $A^J_I:= (A^J)_I$.

If $I \cap J$ is finite, then $A_I^J$ is a finite dimensional quasi-hereditary algebra, 
whenever $\mr A$ is a highest weight category (\cite{CPS:88}, Theorem 3.5).

\begin{prop}\label{updown}
$\mr A$ is a highest weight category if and only if $A^J_I$ is quasi-hereditary for all finitely generated coideals $I$ and finitely generated ideals $J$ such that $I \cap J$ is finite.
\end{prop}

\emph{Proof.}  As noted above, the ``only if'' statement is well known \cite{CPS:88}.  
So suppose $A^J_I$ is quasi-hereditary
for all suitable $I$ and $J$.                                                                        
By a standard argument of Dlab \cite{Dlab}, 
the existence of a highest weight structure on $\mr A$ is equivalent to the surjective multiplication map
$$\frac{Ae_{\lambda}}{\sum_{\mu> \la}Ae_{\mu}
Ae_{\la}}\otimes_F \frac{e_{\la}A}{\sum_{\mu > \la}e_{\la}Ae_{\mu}A}\;\longrightarrow\;\frac
{\sum_{\mu \geq \la}Ae_{\mu}A}{\sum_{\mu > \la}Ae_{\mu}A}$$
being an isomorphism, for all $\lambda \in \Lambda$.
But this can be checked on arbitrarily large finite truncations of $\Lambda$ containing $\lambda$
$\square$

\begin{cor}\label{leftright}
For a locally finite-dimensional algebra $A$, $A \ml$  is a highest weight category if and only if $\mr A$ is a highest weight category.
\end{cor}

\emph{Proof.} 
Follows immediately from Proposition \ref{updown} and the same statement for finite-dimensional algebras
(\cite{ParshallScott}, 4.3(b)).$\square$

\begin{defn}
A locally finite-dimensional algebra $A$ is \emph{quasi-hereditary} 
if $A \ml$ and $\mr A$ are highest weight categories.
\end{defn}

Note that by corollary~\ref{leftright}, 
we can now move freely between left and right modules, standard and  costandard filtrations and we have the usual duality relations between standard modules on one side and costandard modules on the other: 
$\Delta^r(\lambda) \cong \nabla(\lambda)^*, \nabla^r(\lambda) \cong \Delta(\lambda)^*$.

For the rest of this chapter, let $A$ be a locally finite-dimensional quasi-hereditary algebra with poset $\Lambda$ of weights, left standard modules $\Delta(\la)$, left costandard modules $\nabla(\la)$, right standard modules $\Delta^{r}(\la)$ and  right costandard modules $\nabla^{r}(\la)$.
The remaining propositions in this chapter are all proved by cutting down to a suitable finite-dimensional subquotient and applying Ringel's tilting theory for 
finite-dimensional quasihereditary algebras there \cite{Ringel}. We therefore omit the proofs. 

\begin{defn}
$T \in A \ml$ is called \emph{tilting} if it is filtered by standard and by costandard modules.
\end{defn}

\begin{prop}
There is a one-to-one correspondence between $\Lambda$ and the set of indecomposable tilting modules in $A \ml$.
\end{prop}

We denote by $T(\lambda)$ the unique indecomposable tilting module such that $[T(\lambda): L(\lambda)] = 1$,
and $[T(\lambda): L(\mu)] \neq 0$ implies $\mu \leq \lambda$. 

\begin{defn}
We say that $A'$ is \emph{Ringel dual} to $A$ if there exist multiplicities $n_\lambda \in \mathbb{Z}_{\geq 1}$, 
such that $A' \cong \underset{\la, \mu \in \Lambda}{\bigoplus} \Hom_A(T(\la)^{n_\lambda},T(\mu)^{n_\mu})$.
\end{defn}

If $A, A'$ are Ringel dual, then $T = \bigoplus_{\lambda \in \Lambda} T(\la)^{n_\la}$ is an $A$-$A'$ bimodule.
In these circumstances, we call it a tilting bimodule.

For any subset $\Gamma$ of $\Lambda$ let $\Gamma'$ equal $\Gamma$ as a
set, but with the opposite order. Thus, for an ideal $J \subseteq
\Lambda$ we obtain a coideal $J' \subseteq \Lambda' $, for a coideal
$I \subseteq \Lambda$ we obtain an ideal $I' \subseteq \Lambda' $.

\begin{prop}
$A'$ is quasi-hereditary with poset $\Lambda'$.
\end{prop}

\begin{prop}\label{doubledual} $A'' \cong A$.
\end{prop}

\begin{prop}\label{cutdual}

${}$

(i) $(A^J)' \cong  A'_{J'}$ 

(ii)$(A_I)' \cong (A')^{I'}$.
\end{prop}

\section{Algebraic constructions}

Throughout this chapter $A$ will be a finite-dimensional algebra,
endowed with an $A$-$A$-bimodule $T$.

Define $B_0 := \underset{i\in \ZZ}{\bigoplus} A_i$ where $A_i \cong A$
for all $i \in \ZZ$. 
We define $B_1 :=\underset{i\in \ZZ}{\bigoplus} {}_i T_{i+1}$ as a $B_0,B_0$-bimodule, 
where each ${}_i T_{i+1}$ is isomorphic to $T$ but with action of $A_i$ on the left and of $A_{i+1}$ on the right. 

Let $B$ be the trivial extension of $B_0$ by $B_1$; we can think of this as a matrix

$$B= \left(   
\begin{array}{cccccc}
\ddots & {}_{i-2}T_{i-1} & 0 &  \cdots & \\
0& A_{i-1} & {}_{i-1}T_i& 0 & \cdots&\\
\cdots & 0 &  A_{i} & {}_{i}T_{i+1}& 0& \cdots\\
 &\cdots  &0&  A_{i+1} & {}_{i+1}T_{i+2}& 0\\
 & &&& A_{i+2}  &\ddots\\
&& &&&  \ddots \\
\end{array}
\right)$$
where the $A_i$ are on the leading diagonal.
Let $$B^* = \bigoplus_{i \in \mathbb{Z}} \Hom_F(B1_{A_i},F),$$
a $B$-$B$-bimodule. 
Let $C$ be the trivial extension of $B$ by $B^*$. 
Then $C$ is a locally finite dimensional, symmetric algebra.

Let $C^n$ denote the quotient $C/ \sum_{k>n} C1_{A_i}C$ of $C$.
Let $C_1^n$ denote the subalgebra $\sum_{i,j \geq 1} 1_{A_i} C^n 1_{A_j}$ 
of $C^n$. 

\begin{lem} \label{Cstructure}
The algebra $C_1^n$ is $\mathbb{Z}$-graded, concentrated in degrees $0$,$1$, and $2$.
In descending vertical order, its components in degrees $0$, $1$ and $2$ are,
$$\bigoplus_{1 \leq i \leq n} A_i$$                         
$$\bigoplus_{1 \leq i \leq n-1} ({_iT_{i+1}} \oplus {_iT_{i+1}}^*)$$                             
$$\bigoplus_{1 \leq i \leq n-1} A_i^*.$$
\end{lem}

\begin{lem}
Suppose that $T \cong T^*$, as $A$-$A$-bimodules.
Then the infinite dihedral group $D_\infty$ acts as automorphisms of $C$.
The space 
$$T_1^n = \underset{\substack{0 \leq j \leq n-1 \\1 \leq i \leq n }}{\bigoplus}1_{A_i}C1_{A_j}$$ 
has the structure of a self-dual $C_1^n$-$C_1^n$-bimodule.
\end{lem}
\emph{Proof.}

We define an action of
$D_\infty$ on $C$ as follows: The involution $\sigma$ sends
${A_i}$ to ${A_{-i}}$ via the identity, \;
${A^*_i}$ to ${A^*_{-i}}$ via the identity, \;
${}_{i-1}T_i$ to ${}_{-(i-1)} T^*_{-i}$ via the isomorphism $T \cong T^*$, 
and analogously ${}_{i+1}T^*_{i}$ to ${}_{-(i+1)}T_{-i}$.  
Thanks to the assumption that ${}_A T^*_A \cong {}_A T_A$, we see that this
is indeed an algebra isomorphism.
The
translation $\tau$ in $D_\infty$ just maps the components 
$C 1_{A_i}$ to the analogous components of $C 1_{A_{i+1}}$, which
is also clearly an isomorphism.

Of course, $C$ itself is a $C$-$C$-bimodule, 
but what about the truncation $T_1^n$ ?
The idempotents $1_{A_i}$, for $i>n$ act as zero on $C1_{A_j}$, for $j<n$.
Therefore, $C_1^n$ acts naturally on the left of $T_1^n$.
After twisting the right action of $C$ on itself 
by the automorphism $\sigma \circ \tau^{-n}$, 
we can similarly observe a right action of $C_1^n$ on $T_1^n$.

Now 
$$(\underset{\substack{0 \leq j \leq n-1 \\1 \leq i \leq n }}
{\bigoplus}1_{A_i}C1_{A_j})^* \cong 
{}(\underset{\substack{0 \leq j \leq n-1 \\1 \leq i \leq n }}
{\bigoplus}1_{A_j}C1_{A_i})$$ 
by self-duality of $C$.
$\square$

\begin{defn}
Let  
$${\cC}_n \circlearrowright \{ \textrm{algebras with a self-dual bimodule} \}$$
be the map which takes the pair $(A,T)$ to the pair $(C_1^n, T_1^n)$.
\end{defn}

When employing the above definition, we sometimes forget the self-dual bimodules, 
and write simply ${\cC}_n(A)$ for the algebra $C_1^n$.

Assume now that $A$ is a quasi-hereditary algebra with poset $\Lambda$ of weights. 
Let $\Lambda_{B}^1 = \amalg_{i\in \ZZ}
\Lambda[i]$ of weights, with the same ordering as in $\Lambda$ within each $\Lambda[i]$, and 
$\la[i] < \mu[j]$ for $i \neq j$ if and only if 
$i > j \in \ZZ$.
Let $\Lambda_{B}^2 = \amalg_{i\in \ZZ}
\Lambda[i]$ of weights, with the same ordering as in $\Lambda$ within each $\Lambda[i]$, and 
$\la[i] < \mu[j]$ for $i \neq j$ if and only if 
$i < j \in \ZZ$.

The partially ordered sets $\Lambda_B^1,\Lambda_B^2$ index the irreducible
$B_0$-modules.
Indeed, $B_0$ is obviously locally finite-dimensional and quasi-hereditary with 
respect to the posets $\Lambda_B^1, \Lambda_B^2$.

For the remainder of this chapter, 
we assume that $A$ is Ringel self-dual, and that $T$ is a tilting bimodule for $A$, such that $T_A \cong (_AT)^*$. 
Thus, $T_A \in {\F}(\Delta^r) \cap {\F}(\nabla^r)$.

\begin{thm}\label{B}
$B$ is quasi-hereditary with respect to the poset $\Lambda_B^1$,
with standard and costandard modules
$$\Delta^1_B(\la[i]) = \Delta_{B_0}(\la[i]) \qquad \hbox{and} \qquad \nabla^1_B(\la[i]) = \Hom_{B_0}(B,\nabla_{B_0}(\la[i])).$$
$B$ is quasi-hereditary with respect to $\Lambda_B^2$,
with standard and costandard modules 
$$\Delta^2_B(\la[i]) = B \underset{B_0}{\otimes}\Delta_{B_0}(\la[i])   \qquad \hbox{and} \qquad \nabla^2_B(\la[i]) = \nabla_{B_0}(\la[i]).$$
$B$ is Ringel self-dual and Ringel duality exchanges the two 
quasi-hereditary structures on $B$.
\end{thm}

\emph{Proof.}

First observe that indeed $B$ is locally finite-dimensional and $\Lambda_B$ 
indexes simple modules since $B_1$ forms a nilpotent ideal in $B$. 

(1) $\Delta_{B}^1(\la[i])$ has a simple top and the radical consists of composition factors with smaller indices.

Obvious from $B_0$.

(2) $B \projl \subset \F(\Delta^1_B)$ with order relations as required.

We show that ${}_B B 1_{A_i} \in \F(\Delta^1_B)$ for all $i \in \ZZ$.
But ${}_B B 1_{A_i}$ has a filtration with a submodule ${}_{i-1}T_i$
as submodule and $A_i$ as quotient. As left $B$-module, ${}_{i-1}T_i$
is filtered by $\Delta_{B_0}(\la[i-1])$ and $A_i$ is filtered by
$\Delta_{B_0}(\la[i])$ with $\la \in \Lambda$. Since for $A_i$ the
filtration by $\Delta_{B_0}$'s is in the right order (on every direct
summand) and $\Delta_{B_0}(\la[i-1]) >\Delta_{B_0}(\mu[i]) $ for all
$\la, \mu \in \Lambda$, the filtration respects the necessary
inequalities on labels.

(3) $\Delta_{B}^2(\la[i])$ has a simple top and the radical consists of composition factors with smaller indices.

$\Delta_{B}^2(\la[i])$ has a submodule isomorphic to $B_1 \underset{B_0}{\otimes}\Delta_{B_0}(\la[i]) \cong {}_{i-1}T_i \underset{A_i}{\otimes} \Delta_{A_i}(\la)$,
the quotient by which is isomorphic to
$B_0 \underset{B_0}{\otimes}\Delta_{B_0}(\la[i]) \cong \Delta_{B_0}(\la[i])$. 
The latter has simple head,
and all other compositions factors have smaller indices by the 
quasihereditary structure of $B_0$. 
The former has composition factors with labels in $\Lambda[i-1]$ which, 
since in this ordering $i-1 <i$, are smaller as desired. 
Furthermore $B_1$ is a nilpotent ideal in $B$, thus the above submodule does not 
contribute to the head of the module.

(4)  $B \projl \subset \F(\Delta^2_B)$ with order relations as required.

$B_{B_0} \cong (B_0)_{B_0} \oplus (B_1)_{B_0}$ and $(B_0)_{B_0}$ is projective, hence flat. We claim that $(B_1)_{B_0} \underset{B_0}{\otimes} -$ is exact on 
$\F(\Delta_{B_0})$. To prove this, it suffices to check that ${}_{i-1}T_i \underset{A_i}{\otimes} -$ is exact on $\F(\Delta_{A_i})$. So let $M \in \F(\Delta_{A_i})$ and consider ${}_{i-1}T_i \underset{A_i}{\otimes} M$. This being in $\F(\Delta_{A_{i-1}})$, is equivalent to $({}_{i-1}T_i \underset{A_i}{\otimes} M)^*$ being in $\F(\nabla^{r}_{A_{i-1}})$. Now
$({}_{i-1}T_i \underset{A_i}{\otimes} M)^* = \Hom_F({}_{i-1}T_i \underset{A_i}{\otimes} M, F) \cong \Hom_{\mr A_i}(T_i, M^*)$.
But $M \in \F(\Delta_{A_i})$ implies $ M^* \in \F(\nabla^{r}_{A_{i}})$ and, 
by the assumption that $T_A \cong ({}_A T)^*$, $T_A$ is also a tilting module for $\mr A$. 
Therefore, $\Hom_{\mr A_i}(T_i,-)$ is exact on $\F(\nabla^{r}_{A_{i}})$ by \cite{D:98}, A4 (1),
and thus $({}_{i-1}T_i \underset{A_i}{\otimes} M)^* \in \F(\nabla^{r}_{A_{i-1}})$.
So $B \underset{B_0}{\otimes} -$ is exact on $\F(\Delta_{B_0})$, 
and ${}_BB \in \F(\Delta_{B}^2)$. 
The required ordering conditions follow immediately from those for $B_0$.

This finishes the proof of $B$ having two quasihereditary structures.

Similarly, we find that for the right module categories, with respect to $\Lambda_B^1$, we have
$$\Delta^{1,r}_B(\la[i]) =  \Delta^{r}_{B_0}(\la[i])\underset{B_0}{\otimes} B$$ 
and with respect to $\Lambda_B^2$,
$$\Delta^{2,r}_B(\la[i]) = \Delta^{r}_{B_0}(\la[i]) $$

By duality, we now see that 
\begin{equation*}\begin{split}\nabla^1_B(\la[i]) &= (\Delta^{1,r}(\la[i]))^* = \Hom_{F}(\Delta^{r}_{B_0}(\la[i])\underset{B_0}{\otimes} B,F) \\ &\cong \Hom_{B_0}(B, (\Delta^{r}_{B_0}(\la[i]))^*) \cong  \Hom_{B_0}(B, \nabla_{B_0}(\la[i])),\end{split}\end{equation*}
$$\nabla^{1,r}(\la[i]) = (\Delta^{1}(\la[i]))^* = (\Delta_{B_0}(\la[i]))^* = \nabla^{r}_{B_0}(\la[i]),$$
$$\nabla^2_B(\la[i]) =(\Delta^{2,r}(\la[i]))^* = (\Delta^{r}_{B_0}(\la[i]))^* =  \nabla_{B_0}(\la[i]),$$
and
\begin{equation*}\begin{split}\nabla^{2,r}(\la[i])&= (\Delta^{2}(\la[i]))^* = \Hom_{F}( B\underset{B_0}{\otimes}\Delta_{B_0}(\la[i]),F) \\&\cong \Hom_{\mr B_0}(B, (\Delta_{B_0}(\la[i]))^*) \cong  \Hom_{B_0}(B, \nabla^{r}_{B_0}(\la[i])).\end{split}\end{equation*}

To prove the Ringel self-duality of $B$, we need the following lemma.

\begin{lem}
$\Delta^2_B(\la[i]) \cong \nabla^1_B(\la'[i-1]).$
\end{lem}

\emph{Proof of the lemma.} We know that 
$(\nabla^1_B(\la'[i-1]))^* \cong \Delta^{1,r}_B(\la'[i-1]))$, 
so it suffices to show that there exists a non-degenerate bilinear form 
$\langle \;,\; \rangle: \Delta^{1,r}_B(\la'[i-1])) \times \Delta^2_B(\la[i]) \longrightarrow F$ 
with the property that 
$\langle x,by \rangle= \langle xb,y \rangle$ 
for $x \in \Delta^{1,r}_B(\la'[i-1])), y \in \Delta^2_B(\la[i]), b \in B$. 
This is equivalent to having a linear map 
$\Delta^{1,r}_B(\la'[i-1])) {}_F \otimes \Delta^2_B(\la[i]) \longrightarrow F$
which factors over
\begin{equation*}\begin{split}
\Delta^{1,r}_B(\la'[i-1])) \otimes_B \Delta^2_B(\la[i]) 
&= \Delta^{r}_{B_0}(\la'[i-1])\underset{B_0}{\otimes} B \underset{B}{\otimes} B \underset{B_0}{\otimes}\Delta_{B_0}(\la[i])\\ 
& \cong \Delta^{r}_{B_0}(\la'[i-1])\underset{B_0}{\otimes} B \underset{B_0}{\otimes}\Delta_{B_0}(\la[i])\\
& \cong \Delta^{r}_{B_0}(\la'[i-1])\underset{B_0}{\otimes} B_0 \underset{B_0}{\otimes}\Delta_{B_0}(\la[i])\\& \quad \oplus \Delta^{r}_{B_0}(\la'[i-1])\underset{B_0}{\otimes} B_1 \underset{B_0}{\otimes}\Delta_{B_0}(\la[i]).
\end{split}\end{equation*}

But 
$\Delta^{r}_{B_0}(\la'[i-1])\underset{B_0}{\otimes} B_0 \underset{B_0}{\otimes}\Delta_{B_0}(\la[i]) \cong \Delta^{r}_{B_0}(\la'[i-1])1_{A_{i-1}}\underset{B_0}{\otimes}1_{A_{i}}\Delta_{B_0}(\la[i]) \cong 0$ 
and we claim that
$\Delta^{r}_{B_0}(\la'[i-1])\underset{B_0}{\otimes} B_1 \underset{B_0}{\otimes}\Delta_{B_0}(\la[i])$ 
is isomorphic to 
$ \nabla^{r}_{B_0}(\la[i])\underset{B_0}{\otimes}\Delta_{B_0}(\la[i])$ 
on the one hand and to 
$\Delta^{r}_{B_0}(\la'[i-1])\underset{B_0}{\otimes} \nabla_{B_0}(\la'[i-1])$ 
on the other hand.
Indeed, 
\begin{equation*}\begin{split}
({}_{i-1}T_i \underset{A_i}{\otimes}\Delta_{A_i}(\la))^* 
& = \Hom_F( {}_{i-1}T_i \underset{A_i}{\otimes}\Delta_{A_i}(\la) ,F) 
\cong \Hom_{\mr A_i}({}_{i-1}T_i , (\Delta_{A_i}(\la))^*  )\\ 
& \cong \Hom_{\mr A_i}({}_{i-1}T_i, \nabla^{r}_{A_i}(\la))  
\cong \Delta^{r}_{A_{i-1}}(\la'),
\end{split}\end{equation*}
whence 
${}_{i-1}T_i \underset{B_0}{\otimes}\Delta_{B_0}(\la[i]) ={}_{i-1}T_i \underset{A_i}{\otimes}\Delta_{A_i}(\la) \cong \nabla_{A_i}(\la') = \nabla_{B_0}(\la'[i-1])$. The second isomorphism,
$\Delta^{r}_{B_0}(\la'[i-1])\underset{B_0}{\otimes} B_1 \cong \nabla^{r}_{B_0}(\la[i]) $ 
follows by the same arguments for right modules.  But
$\nabla^{r}_{B_0}(\la[i])$ and $\Delta_{B_0}(\la[i])$ are dual to one
another as are $\Delta^{r}_{B_0}(\la'[i-1])$ and
$\nabla_{B_0}(\la'[i-1])$. Thus we have unique (up to scalar)
non-degenerate bilinear forms on both pairs.  Defining our bilinear
form as the sum of both gives us a $B$-equivariant nondegenerate
bilinear form on $\Delta^{1,r}_B(\la'[i-1])) \times
\Delta^2_B(\la[i]) $. (It is nondegenerate because, as $B_0$-modules,
$\Delta^{1,r}_B(\la'[i-1])) \cong \Delta^{r}_{B_0}(\la'[i-1]) \oplus \nabla^{r}_{B_0}(\la[i]) $
and 
$\Delta^2_B(\la[i]) \cong \Delta_{B_0}(\la[i]) \oplus \nabla_{B_0}(\la'[i-1])$
and it is nondegenerate on both dual pairs.)
This completes the proof of the lemma. $\square$

\bigskip

\emph{Proof of theorem, continued.}  By the lemma 
$B \projl \subset \F(\nabla^1_B)=\F(\Delta^2_B)$, 
but we also have 
$B \projl \subset \F(\Delta^1_B)$, 
hence projective modules are tilting modules in the first highest weight structure on $B \ml$. But clearly
$$\underset{\la[i], \mu[j] \in \Lambda_B}{\bigoplus} \Hom_B(P(\la[i]),P(\mu[j]))
\cong B$$
so $B$ is indeed Ringel self-dual.
Denoting the new standard modules by $\tilde{\Delta}_B$ we obtain
$$\tilde{\Delta}^1_B(\la[i]) = \Hom_B(B, \nabla^1_B(\la[i])) \cong \nabla^1_B(\la[i]) \cong \Delta^2_B(\la'[i+1]).$$
By the right analogue of the lemma we see that (right) projectives are tilting modules for the second highest weight structure on $\mr B$, and by the same computation as above, we obtain
$$\tilde{\Delta}^{2,r}_B(\la[i]) = \Hom_B(B, \nabla^{2,r}_B(\la[i])) \cong \nabla^{2,r}_B(\la[i]) \cong \Delta^{1,r}_B(\la'[i-1]).$$
Dualizing we see that 
$$\tilde{\nabla}^{1,r}_B(\la[i]) \cong (\tilde{\Delta}^1_B(\la[i]))^* 
\cong (\Delta^2_B(\la'[i+1]))^* \cong \nabla^{2,r}_B(\la'[i+1])$$
and
$$\tilde{\nabla}^{2}_B(\la[i]) \cong (\tilde{\Delta}^{2,r}_B(\la[i]))^*
\cong (\Delta^{1,r}_B(\la'[i-1]))^* \cong \nabla^1_B(\la'[i-1]).$$
Since $\Delta$'s and $\nabla$'s determine each other, this completes the proof of the theorem.
$\square$

\bigskip

Let $\Lambda_C^1 = \Lambda_B^1$, and $\Lambda_C^2 = \Lambda_B^2$.

\begin{thm}\label{C}
$C$ is quasi-hereditary with poset $\Lambda_C^1$, as well as with 
poset $\Lambda_C^2$. We have
$$\nabla_C^1(\la[i]) = \nabla_B^1(\la[i]), \qquad \Delta_C^2(\la[i]) = \Delta_B^2(\la[i]).$$
Furthermore, $C$ is Ringel self-dual, and Ringel duality exchanges the two highest weight 
structures on $C$.
\end{thm}

\emph{Proof.}  The equality of the indexing sets for simple modules follows from the
nilpotency of $B^*$ in $C$.  Now, $C 1_{A_i}$ has a filtration with
submodule $B^*1_{A_i} \cong (1_{A_i}B)^*$ and quotient $B1_{A_i}$. The latter has a
filtration by $\Delta_B^2(\la[i])$, where $\la \in \Lambda$, with the
necessary properties by Theorem \ref{B}.  
The former has a filtration by $(\Delta^{1,op}_B(\la[i]))^* \cong \nabla^1_B(\la[i]) \cong \Delta_B^2(\la'[i+1])$. 
So, since $i+1 >i$ we have a filtration respecting the necessary inequalities on labels.

The fact that $C$ is symmetric follows from a general statement the the trivial extension of an algebra by its dual is symmetric.

Ringel self-duality follows immediately from symmetry, since projectives have a $\Delta$-filtration, but since they are the same as injectives, also a $\nabla$-filtration, thus projectives are tilting modules, implying Ringel self-duality.
$\square$

\bigskip

Set $J_n:= \underset{j \leq n}{\bigcup} \Lambda[j]$ and $I_k:= \underset{i \geq k}{\bigcup} \Lambda[i]$ and
adopt the notational convention $C^{n}:=C^{J_n}$, $C_k:= C_{I_k}$, and $C_k^{n}:=C_{I_k}^{J_n}$.

Let us now assume that ${}_A T^*_A \cong {}_A T_A$ as a bimodule.
Recall that in these circumstances,
$D_\infty= <\sigma,\tau>$ acts on $C$.  
Note that in the Ringel duality in
theorem \ref{C}, $C' = \tau^{-1}(C)$, since the projective
$P(\la[i])$ has a submodule $\Delta_C(\la'[i+1])$, implying
$P(\la[i]) \cong T_C(\la[i+1])$ and $P_{C'}(\la[i]) =  \Hom_C(\underset{\substack{j \in \ZZ \\ \la \in \Lambda}}{\bigoplus} P(\la[j]), P(\la[i-1]))$.

\begin{thm}\label{C[1,n]}
$C_1^n$ for $n \geq 1$ is Ringel self-dual, and the tilting bimodule 
$T_{C_1^n}$ is a self-dual bimodule.
\end{thm}

\emph{Proof.} By Proposition \ref{cutdual}, $(C_1^n)' \cong (C')_{J_n'}^{I_1'} $
with the ordering $i > i+1$ on $\ZZ$.  
Therefore 
\begin{equation*}\begin{split}(C_1^n)' &\cong (\tau^{-1}C)_{J'_n}^{I'_1} \cong C_{J'_{n-1}}^{I'_0}  \overset{\sigma}{\cong}
C^{0}_{-(n-1)} \overset{\tau^{n}}{\cong} C^n_{1}.\end{split}\end{equation*}

The tilting module $T_{C_1^n}$ satisfies
\begin{equation*}\begin{split}T_{C_1^n} & = \underset{\substack{j \leq n \\ \la \in \Lambda}}{\bigoplus} \Hom_C(\underset{\substack{i \geq 1}}{\bigoplus} C1_{A_i} , T_C(\la [j])) \\ 
& \cong \underset{\substack{j \leq n \\ \la \in \Lambda}}{\bigoplus} \Hom_C(\underset{\substack{i \geq 1 }}{\bigoplus} C1_{A_i}  , P_C(\la [j-1]))\\
& \cong \underset{\substack{j \leq n-1 \\ \la \in \Lambda}}{\bigoplus}\Hom_C(\underset{\substack{i \geq 1 }}{\bigoplus} C1_{A_i} ,  P_C(\la [j])) \\
& = \underset{\substack{0 \leq j \leq n-1 \\1 \leq i \leq n }}{\bigoplus}\Hom_C(C1_{A_i} , C1_{A_j} )\\
& \cong \underset{\substack{0 \leq j \leq n-1 \\1 \leq i \leq n }}{\bigoplus}1_{A_i}C1_{A_j}
 .\end{split}\end{equation*}
The first equality comes from the fact that factoring out a heredity ideal doesn't change the tilting modules for the remaining labels and that the tilting module for a heredity subalgebra is the tilting module multiplied by the idempotent.The fourth equality takes into account that we only have nonzero maps from $C1_{A_i}$ to itself or to $C1_{A_{i \pm 1}}$.
Now $$({}_{C_1^n}(\underset{\substack{0 \leq j \leq n-1 \\1 \leq i \leq n }}{\bigoplus}1_{A_i}C1_{A_j})_{(C_1^n)'})^* = {}_{(C_1^n)'}(\underset{\substack{0 \leq j \leq n-1 \\1 \leq i \leq n }}{\bigoplus}1_{A_j}C1_{A_i})_{C_1^n}$$ (by self-duality of $C$), but to view this as a $(C_1^n, (C_1^n)')$-bimodule we have to twist with $\sigma \circ \tau^{-n}$ on the left and its inverse on the right which yields $\underset{\substack{0 \leq j \leq n-1 \\1 \leq i \leq n }}{\bigoplus}1_{A_i}C1_{A_j}$ as desired.
$\square$

\begin{cor} The map ${\cC}_n$ restricts to a map
$${\cC}_n \circlearrowright 
\{ \textrm{quasi-hereditary algebras with a self-dual tilting bimodule} \}. 
\quad \square$$
\end{cor}

\section{Schur algebras} \label{Schur}

Let $M$ denote the algebra of $n \times n$ matrices over $F$.
Recall the Schur algebra ${\cS}(n,r)$ is defined to be the subalgebra 
$(M^{\otimes r})^{\Sigma_r}$ of fixed points under the action of the symmetric group 
$\Sigma_r$ on $M^{\otimes r}$.
The category of representations of ${\cS}(n,r)$ can be identified with the 
category of polynomial representations of $GL_n(F)$, of degree $r$ \cite{Green}.

Let $\Lambda^+(n,r)$, the set of partitions of $r$ with $n$ parts or fewer, given the dominance ordering.
The algebra ${\cS}(n,r)$ is quasi-hereditary with respect to the poset $\Lambda^+(n,r)$.
We write $\xi_\lambda \in {\cS}(n,r)$ for Green's idempotents in ${\cS}(n,r)$, 
for $\lambda \in \Lambda(n,r)$.

In this paper, we are only concerned with ${\cS}(2,r)$, 
but it will be useful to recall some facts about Ringel duality which hold for general $n$.

\begin{lem} (S. Donkin, \cite{D:98}, 4.1)
Let $n \geq r$.
Then $\bigwedge^r(M)$ is a tilting ${\cS}(n,r)$-${\cS}(n,r)$-bimodule. $\square$
\end{lem}

When $n \geq r$, let ${\cS}'(n,r) = {\cS}(n,r)$. 
When $n<r$, 
let $${\cS}'(n,r) \cong {\cS}(r,r)/ \sum_{\lambda \notin \Lambda(n,r)} {\cS}(r,r) \xi_\lambda {\cS}(r,r).$$
The algebras ${\cS}(n,r), {\cS}'(n,r)$ are Ringel dual.

The Schur algebra possesses a natural anti-automorphism inherited from the transpose operator on $M$.
We also call this antiautomorphism the transpose operator, and denote by $s^T$ the twist of an element $s$ 
by the transpose operator.
Since $\xi_\lambda^T = \xi_\lambda$ for all $\lambda$, 
the transpose operator descends to an antiautomorphism of ${\cS}'(n,r)$

If $A$ is an algebra, endowed with an antiautomorphism $x$, then given any left/right $A$-module $M$, 
we define the right/left $A$-module $M^{op}$ to be that obtained by twisting the action of $A$ on $M$ by $x$.
If $A_1,A_2$ are algebras, endowed with antiautomorphisms $x_1,x_2$, then given an $A_1$-$A_2$-bimodule $M$,
we define the $A_2$-$A_1$-bimodule $M^{op}$ to be that obtained by twisting the actions of 
$A_1,A_2$ on $M$ by $x_1,x_2$.

\begin{lem} \label{tiltdual}
Let $_{{\cS}(n,r)}T_{{\cS}'(n,r)}$ be a tilting bimodule. Then $T^{op} \cong T^*$, as 
${\cS}'(n,r)$-${\cS}(n,r)$-bimodules, 
where $T^{op}$ is obtained after twisting 
$T$ by the transpose antiautomorphisms of ${\cS}'(n,r)$, ${\cS}(n,r)$. 
\end{lem}

\emph{Proof.} In case $n \geq r$, we have $T \cong \bigwedge^r(M)$.
However, it is well known that $\bigwedge^r(M)$ is self-dual, 
which is to say $\bigwedge^r(M)^{op} \cong \bigwedge^r(M)^*$.

The case $n<r$ follows by truncation from the case $n=r$.
Indeed, in this case, we have $T = (\sum_{\lambda \in \Lambda(n,r)} \xi_\lambda) \bigwedge^r(M)$.
Since $\xi_\lambda^T = \xi_\lambda$, this bimodule is also self-dual. $\square$

\bigskip

We now restrict our study to the case $n=2$.
Suppose $F$ is a field of characteristic $p>0$.
Let $\cS= \cS(2,r)$ be the Schur algebra over $F$, 
where $r= ap^k -2$ or $r= ap^k -3$ for some $k \geq 1$ and $2 \leq a \leq p$. 
Along with the cases $r< p^2$ and $r= ap^k -1$, these are exactly the Schur algebras which are Ringel self-dual
(\cite{EH:02}, Theorem 27).
Furthermore, 
$\cS(2,ap^k - 1)$ is Morita equivalent to $\cS(2,ap^k - 3) \oplus F$ (\cite{EH:02}, Corollary 2).

If $r$ is odd, our index set $\Lambda$ for the quasi-hereditary structure of $\cS$ 
consists of all odd natural numbers up to $r$; 
if $r$ is even, it consists of all even natural numbers up to $r$, including $0$.

\bigskip

The following definitions assume $p$ odd.
If $r$ is odd, let $A= \cS(2, p^k-2)$ and if $r$ is even, let $A= \cS(2, p^k-3)$.
We define subsets $I_j$, for $1 \leq j \leq a$, of $\Lambda$ as follows:

\begin{center}
\begin{tabular}{c|c}
         &  $r$ odd               \\ \hline
$j$ odd  & \small{$I_j = \{\la \in \Lambda \mid (j-1)p^k +1 \leq \la \leq jp^k-2 \}$} 
\\ \hline
$j$ even & \small{$I_j = \{\la \in \Lambda \mid (j-1)p^k \leq \la \leq jp^k-3 \}$}
         
\\ \hline
         &  $r$ even              \\ \hline
$j$ odd  &  \small{$I_j = \{\la \in \Lambda \mid (j-1)p^k \leq \la \leq jp^k-3 \}$}
\\ \hline
$j$ even &  \small{$I_j = \{\la \in \Lambda \mid (j-1)p^k +1 \leq \la \leq jp^k-2 \}$}
\end{tabular}
\end{center}
\bigskip

In case $p=2$ (and thus necessarily $a=2$), let $A = \cS(2,2^k-3)$ if $r$ is odd, 
and $A = \cS(2, 2^k-2)$ if $r$ is even.
We define subsets $I_j$, for $j=1,2$, of $\Lambda$ as follows:

\begin{center}
\begin{tabular}{c|c}
         &  $r$ odd               \\ \hline
$j=1$  & \small{$I_j = \{ 1,3,\dots,2^{k-1}-3  \}$} 
\\ \hline
$j=2$ & \small{$I_j = \{\ 2^{k-1}+1,2^{k-1}+3, \dots,2^k-3 \}$}
         
\\ \hline
         &  $r$ even              \\ \hline
$j=1$  &  \small{$I_j = \{ 0,2,\dots,2^{k-1}-2 \}$}
\\ \hline
$j=2$ &  \small{$I_j = \{ 2^{k-1},2^{k-1}+2, \dots,2^k-2 \}$}
\end{tabular}
\end{center}

\bigskip

Let us define $I_0 := \Lambda \setminus (\underset{1 \leq j \leq a}{\bigcup}I_j)$.


For $1 \leq j \leq a$, set $b_j :=\min\{I_j\}$, $r_j := \max\{I_j\}$

\bigskip

We choose orthogonal idempotents $\{ e_\lambda \}_{\lambda \in \Lambda}$ in $\cS$, such that 
$\cS \cong \bigoplus_{\lambda,\mu \in \Lambda} e_\lambda {\cS} e_\mu$, 
and $\cS/{\mathcal J}(\cS) = \bigoplus_{\lambda \in \Lambda} M_\lambda$ 
is a direct sum of matrix rings $M_\lambda$ over $F$, where $e_\lambda$ is the unit of $M_\lambda$.

Let $f_j := \sum_{\la \in I_j} e_\la$, 
where $e_\la \in \cS$ is the primitive idempotent corresponding to the projective $P(\la)$. 
Let $\eps_j = \sum_{i \geq j } f_i$.

\bigskip

By work of A. Henke and S. Koenig,
there are idempotents $\eta_j \in {\cS}$ (denoted $\xi_l^o$ in \cite{HK:02}),
and explicit isomorphisms $\Phi_j: A \rightarrow \eta_{j} \cS \eta_{j} / \eta_{j} \cS \eta_{j+1} \cS \eta_{j}$,
for $1 \leq j \leq a$ (\cite{HK:02}, Theorem 3.3).

We now assume that the idempotents $e_\lambda$ are chosen in such a way that $e_\lambda \eta_j = \eta_j e_\lambda$,
for $\lambda \in \Lambda, 1 \leq j \leq a$. 
It therefore follows that the idempotents $\eta_j$ commute with $f_i, \epsilon_i$ as well,
and we have $\epsilon_k \eta_j = \eta_k$, for $1 \leq j \leq k \leq a$

Let $\alpha_j := f_j \cS f_j / f_j \cS f_{j+1} \cS f_j$, for $1 \leq j \leq a$.

\begin{lem}\label{remark1}
The algebra $\alpha_j$ is Morita equivalent to $A$, for $1 \leq j \leq a$.
We have $f_j \cS f_i =0$ unless $j-1 \leq i \leq j+1$, and
\begin{equation}\begin{split}\label{idempdecomp}
\cS = \bigoplus_{j=1}^a f_j \cS f_j \oplus \bigoplus_{j=1}^{a-1} (f_j \cS f_{j+1}  +f_{j+1} \cS f_j) \oplus 
\bigoplus_{\la \in I_0} e_\la \cS e_\la.
\end{split}\end{equation}
\end{lem}

\emph{Proof:}
From the decomposition matrix of ${\cS}$, 
we see that $f_j \cS f_i =0$ unless $j-1 \leq i \leq j+1$ 
and that for $\la \in I_0$, $e_\la \cS e_\mu =  e_\mu \cS e_\la =0$ unless $\mu = \la$ when it is isomorphic to $F$.
Hence 
$$\alpha_j = f_j \cS f_j / f_j \cS f_{j+1} \cS f_j \cong 
\eps_{j} \cS \eps_{j} / \eps_{j} \cS \eps_{{j+1}} \cS \eps_{j},$$ 
This algebra is Morita equivalent to  
$\eta_{j} \cS \eta_{j} / \eta_{j} \cS \eta_{j+1} \cS \eta_{j}$,
which is isomorphic to $A$. 
This completes the proof of the lemma.
$\square$

\begin{remark}
It will be important to us that the Henke-Koenig isomorphism $\Phi_j$ between 
$A$ and $\eta_{j} \cS \eta_{j} / \eta_{j} \cS \eta_{j+1} \cS \eta_{j}$ 
is compatible with the transpose operators on $\cS, A$. 
To be more explicit, $\eta_j^T = \eta_j$, and $\Phi_j(a^T) = \Phi_j(a)^T$, for $a \in A$.
\end{remark}

\begin{lem}\label{backforth}
We have $f_{a} \cS f_{a-1} \cS f_{a} =0$.
\end{lem}

\emph{Proof.} This is a reformulation of \cite{EH:02}, Proposition 25. 
Indeed, according to this proposition, $\cS f_{a-1} \cS f_{a}$ is the submodule of $\cS f_{a}$ 
consisting of all composition factors of the form $L(\la), \la \in I_{a-1}$, 
implying $$f_{a} \cS f_{a-1} \cS f_{a} \cong \Hom_\cS (\cS f_{a} ,\cS f_{a-1} \cS f_{a}) = 0. \quad \square$$

\begin{lem}\label{tilt}
$X_j = f_j \cS f_{j+1}$ is an $\alpha_j$-$\alpha_{j+1}$-tilting bimodule.
\end{lem}

\emph{Proof.} 
By lemmas \ref{remark1} and \ref{backforth}, the $f_j \cS f_j$-$f_{j+1} \cS f_{j+1}$-bimodule $X_j$
is in fact an $\alpha_j$-$\alpha_{j+1}$-bimodule.
It remains to show that $_{\alpha_j}X_j$ is a full tilting module, and $End_{\alpha_j}(X_j) = \alpha_{j+1}$.

By the same argument as in lemma \ref{remark1} we can
reduce to the case where $a=2$ by considering all modules for the
subalgebra $\eps_{j} \cS \eps_{j} /
\eps_{j} \cS \eps_{{j+2}} \cS \eps_{j}$.
So let $\cS= \cS(2,r)$ where $r \in \{2p^k-2, 2p^k-3\}$ and use the
notation from above. We need to show that $f_1\cS f_2 \in \cS/\cS
\eps_{2} \cS \ml$ is a tilting module. But by
\cite{EH:02}, Proposition 25, $\cS f_1\cS f_2 \subseteq \cS f_2$ is the
submodule consisting of all composition factors of the form $L(\la)$
for $\la \in I_1$ and is isomorphic to the full tilting module for
$\cS(2,\max\{I_1\})$. But by the first of these facts the action
factors over $\alpha_1 = \cS/\cS \eps_{2} \cS \cong 
{\cS}(2,r_1)$, so it is a full tilting module for this algebra.

Now we have a canonical map from $\alpha_2 = f_2 \cS f_2$ to
$\End_{\alpha_1}(f_1\cS f_2)$. Given the fact that $A$ is Ringel self-dual, 
we know that $\alpha_2$, $A$, and $\End_{\alpha_1}(f_1\cS f_2)$ are Morita equivalent, 
thus $\alpha_2$ and 
$\End_{\alpha_1}(f_1 \cS f_2)$ are isomorphic.
It therefore suffices to prove
injectivity of this map. So, suppose it has a nontrivial kernel. This
is equivalent to the existence of an endomorphism $\phi$ of $\cS f_2$,
annihilating all composition factors of the form $L(\la)$ for $\la \in
I_1$ (namely $\cS f_1\cS f_2$). But all composition factors of the
socle of $\cS f_2$ are of the form $L(\la)$
for $\la \in I_1$, by \cite{EH:02}, Lemma 3,
and thus $\im \phi \cap \soc \cS f_2 =0$ forcing $\phi$ to be zero. $\square$

\begin{remark}\label{contained}
Note that it follows from the proof of the lemma that 
$f_j \cS f_{j+1} \cS f_j$ is the annihilator of $f_j \cS f_{j+1}$ 
in $f_j \cS f_j$. Since by Remark \ref{remark1} 
$f_{j-1} \cS f_j \cS f_{j+1} \subseteq f_{j-1} \cS f_{j+1} =0$, 
it follows that $f_j \cS f_{j-1} \cS f_j \subseteq f_j \cS f_{j+1} \cS f_j$.
\end{remark}

Let $\bar{X_j} = f_{j+1} \cS f_j$.
By lemmas \ref{remark1} and \ref{backforth}, $\bar{X_j}$ is an $\alpha_{j+1}$-$\alpha_j$-bimodule.

Let $X_j^{op}$ be the $\alpha_{j+1}$-$\alpha_j$-bimodule 
obtained by passing $_{\alpha_j}X_{j \alpha_{j+1}}$ via the established Morita equivalences to the category of 
$A$-$A$-bimodules, twisting on both sides by the transpose automorphism of $A$, 
and then passing via Morita equivalence to the category of $\alpha_{j+1}$-$\alpha_{j}$-bimodules.

\begin{lem} \label{tiltop}
There is an isomorphism of $\alpha_{j+1}$-$\alpha_j$-bimodules, $\bar{X_j} \cong X_j^{op}$.
\end{lem}

\emph{Proof.} 

We have 
$$X_j = f_j {\cS} f_{j+1} \cong 
\epsilon_j {\S} \epsilon_{j+1}/\epsilon_{j+1} {\cS} \epsilon_{j+1}.$$
This passes, via Morita equivalence, to the $A$-$A$-bimodule
$$\eta_j \cS f_j \underset{f_{j} \cS f_{j}}{\otimes}
f_j \cS f_{j+1} \underset{f_{j+1} \cS f_{j+1}}{\otimes} f_{j+1} \cS \eta_{j+1} \cong$$
$$\eta_j f_j \cS f_{j+1} \eta_{j+1} 
\cong \eta_j {\cS} \eta_{j+1}/\eta_j \epsilon_{j+1} {\cS} \eta_{j+1}.$$
Since twisting by the transpose operator exchanges the irreducible modules
$L(\lambda), L^r(\lambda)$, the projective $\cS$-modules $\cS f_j$ and $\cS f_j^T$ are isomorphic.
We therefore have
$$\bar{X}_j = f_{j+1} {\cS} f_{j} = 
f_{j+1} {\cS} f_{j+1}^T {\cS} f_j^T \cS f_j \cong$$
$$f_{j+1} {\cS} f_{j+1}^T   \underset{\epsilon_{j+1}^T \cS \epsilon_{j+1}^T}{\otimes}
(\epsilon_{j+1}^T {\cS} \epsilon_j^T /\epsilon_{j+1}^T {\cS} \epsilon_{j+1}^T) 
\underset{\epsilon_{j}^T \cS \epsilon_{j}^T}{\otimes}
f_j^T \cS f_j.$$
This passes, via Morita equivalence, to the $A$-$A$-bimodule
$$\eta_{j+1} \cS f_{j+1}^T \underset{f_{j+1}^T \cS f_{j+1}^T}{\otimes}
f_{j+1}^T \cS f_{j}^T \underset{f_{j}^T \cS f_{j}^T}{\otimes} f_{j}^T \cS \eta_{j+1} \cong$$
$$\eta_{j+1} f_{j+1}^T \cS f_{j}^T \eta_{j} \cong 
\eta_{j+1} {\cS} \eta_{j}/\eta_{j+1} {\cS} \epsilon_{j+1}^T \eta_{j}.$$
Applying the transpose anti-automorphism to $\cS$, we exchange the bimodules 
\newline $\eta_j {\cS} \eta_{j+1}/\eta_j \epsilon_{j+1} {\cS} \eta_{j+1}$
and $\eta_{j+1} {\cS} \eta_{j}/\eta_{j+1} {\cS} \epsilon_{j+1}^T \eta_j$, 
the left and right actions being twisted by the transpose operator.
However, the transpose operator is compatible with the Henke-Koenig isomorphisms,
and therefore an equivalent statement is that passing to the opposite exchanges $X_j$ and $\bar{X}_j$.
We therefore have $\bar{X}_j \cong X_j^{op}$, as required.
$\square$

\bigskip

Let us define
$$\mathcal{N} :=\sum_{j=1}^{a-1} (f_j \cS f_{j+1}  +f_{j+1} \cS f_j + f_j \cS f_{j+1} \cS f_j ),$$
$$\mathcal N_2:=\sum_{j=1}^{a-1}f_j \cS f_{j+1} \cS f_j.$$

\begin{prop} \label{filtration}
We have a filtration by of $\cS$ by ideals,
\begin{equation}\begin{split}\label{filt}
\cS \supset \mathcal{N} \supset \mathcal N^2 \supset 0.
\end{split}\end{equation}
Furthermore $\mathcal N^2 = \mathcal N_2$, and $\mathcal N^3 =0$.
We have isomorphisms of  $\cS$-$\cS$-bimodules,
$$\cS/\mathcal N \cong \underset{1 \leq j \leq a}{\bigoplus} \alpha_j \oplus 
\bigoplus_{\lambda \in I_0} e_\lambda {\cS} e_\lambda,$$
$$\mathcal N/\mathcal N^2 \cong \underset{1 \leq j \leq a-1}{\bigoplus}
\left( X_j \oplus X_j^* \right),$$ 
$$\mathcal N^2 \cong \underset{1 \leq j \leq a-1}{\bigoplus} \alpha_j^*.$$ 
\end{prop}

\emph{Proof.} The first statement as well as $\mathcal N^2 = \mathcal N_2$ and $\mathcal N^3 =0$, 
are easily verified using Lemma \ref{remark1}, Lemma \ref{backforth} and Remark \ref{contained}. 
From (\ref{idempdecomp}) we see that 
$$ \cS/\mathcal N \cong \underset{1 \leq j \leq a}{\bigoplus} f_j \cS f_j/ (f_j \cS f_{j+1} \cS f_j) \cong$$
$$\underset{1 \leq j \leq a}{\bigoplus} \alpha_j \oplus \bigoplus_{\lambda \in I_0} e_\lambda \cS e_\lambda,$$
and by lemmas \ref{tiltdual}, \ref{tilt} and \ref{tiltop},  
$$\mathcal N/\mathcal N^2 \cong  
\underset{1 \leq j \leq a-1}{\bigoplus} (f_j \cS f_{j+1}  +f_{j+1} \cS f_j) \cong$$  
$$\underset{1 \leq j \leq a-1}{\bigoplus} ( X_j + X_j^{op}) \cong
\underset{1 \leq j \leq a-1}{\bigoplus} ( X_j + X_j^*).$$

Now all that is left to show is that $f_j \cS f_{j+1} \cS f_j \cong \alpha_j^*$.
To see this, note that by repeatedly applying Remark \ref{remark1} 
\begin{equation*}\begin{split}f_j \cS f_{j+1} \cS f_j &= f_j \cS \eps_{r_{j+1}} \cS f_j \cong f_j \cS \eps_{r_{j+1}} \underset{ \eps_{r_{j+1}}\cS \eps_{r_{j+1}}}{\otimes}   \eps_{r_{j+1}} \cS f_j \\
& \cong f_j \cS f_{j+1} \underset{\alpha_{j+1}}{\otimes} f_{j+1}\cS f_j
\end{split}\end{equation*}
But 
\begin{equation*}\begin{split}
\Hom_F( f_j \cS f_{j+1} \underset{\alpha_{j+1}}{\otimes} f_{j+1}\cS f_j ,F) 
&\cong \Hom_{\mr \alpha_{j+1}}( X_j, X_j^{op*})\\
& \cong  \End_{\mr \alpha_{j+1}}(X_j) \cong \alpha_j, 
\end{split}\end{equation*}
thus $f_j \cS f_{j+1} \cS f_j \cong \alpha_j^*$ as claimed.
$\square$

\bigskip

Let $C_1^a = {\cC}_1^a(A)$ be the algebra obtained by applying the construction ${\cC}_1^a$ of the previous chapter
to the algebra $A$, and its self-dual bimodule $T$.

\begin{thm} \label{grMorita}
The graded algebra $\cS_{gr}$ associated to the filtration $\cS \supset {\cN} \supset {\cN}^2 \supset 0$ 
is Morita equivalent to $C_1^a \oplus F^{\oplus I_0}$.
\end{thm}

\emph{Proof.}
By Proposition \ref{filtration}, we know that $\cS_{gr}$ is Morita eqivalent to $\tilde{C}_1^a \oplus F^{\oplus I_0}$,
where $\tilde{C}_1^a$ is $\mathbb{Z}$-graded, concentrated in degrees $0$,$1$, and $2$.
In descending vertical order, the components of $\tilde{C}_1^a$ in degrees $0$, $1$ and $2$ are,
$$\bigoplus_{1 \leq i \leq a} \tilde{A}_i$$                         
$$\bigoplus_{1 \leq i \leq n-1} ({_i\tilde{T}_{i+1}} \oplus {_i\tilde{T}_{i+1}^*})$$                             
$$\bigoplus_{1 \leq i \leq n-1} \tilde{A}_i^*,$$
where $\tilde{A}_i$ is isomorphic to $A$, and where $_i\tilde{T}_{i+1}$ is a tilting 
$\tilde{A}_i$-$\tilde{A}_{i+1}$-bimodule.
Twisting the isomorphisms $\tilde{A}_i \cong A$ by automorphisms of $A$ if necessary, 
we may assume that $_i\tilde{T}_{i+1} \cong {_AT_A}$.
We proceed to piece together an algebra isomorphism between $\tilde{C}_1^a$ and $C_1^a$ itself. 

We know from the proof of the previous proposition that multiplication  
$f_j \cS f_{j+1} \underset{F}{\otimes} f_{j+1}\cS f_j \twoheadrightarrow  \alpha_j^*$
is surjective, for $1 \leq j \leq a-1$.
Therefore, multiplication $_j\tilde{T}_{j+1} \underset{F}{\otimes} {_j\tilde{T}_{j+1}^*} \twoheadrightarrow 
\tilde{A}_j^*$ is also surjective.

Since we have a canonical isomorphism 
$_j T_{j+1} \underset{A_{j+1}}{\otimes}  {_jT_{j+1}^*} \cong A_{j}^*$, we consequently obtain an isomorphism 
$\tilde{A}_j^* \cong A^*$ of $A$-$A$-bimodules.

We now claim that multiplication
$_j\tilde{T}^*_{j+1} \underset{F}{\otimes} {_j\tilde{T}_{j+1}} \twoheadrightarrow 
\tilde{A}_{j+1}^*$ is also surjective, for $1 \leq j  \leq a-2$. 
Equivalently, we claim that multiplication 
$f_{j+1}\cS f_j \underset{F}{\otimes}f_j \cS f_{j+1} \rightarrow \alpha_{j+1}^*$
is surjective. 
Indeed, this multiplication is inherited from the left module structure on the maximal submodule $M$ of
$\cS f_{j+1}$ whose composition factors $L(\lambda)$ respect $\lambda \in I_{j+1}$. 
The submodule $M$ has a filtration with submodule $\alpha_{j+1}^*$ and
quotient $f_j \cS f_{j+1}$. 
Note that $\cS f_j$ is a tilting module (\cite{EH:02}, Corollary 21, Lemma 24) and therefore self-dual.  
Therefore $M^{op*}$ is the maximal quotient of $\cS f_{j+1}$ all of whose composition factors 
$L(\lambda)$ respect $\lambda \leq I_{j+1}$. 
$M^{op*}$ has a filtration with submodule $f_{j}\cS f_{j+1}$ and quotient isomorphic to $\alpha_{j+1}$.
However, we know the structure of this module precisely.
For instance, the product $f_j \cS f_{j+1} \underset{F}{\otimes} \alpha_{j+1} \rightarrow f_j \cS f_{j+1}$ 
corresponds to the right action of $T \otimes A \rightarrow T$.
Since the product on $M$ is dual to that on $M$, 
the map 
$f_{j+1}\cS f_j \underset{F}{\otimes}f_j \cS f_{j+1} \twoheadrightarrow \tilde{A}_{j+1}^*$
is dual to the map $A \hookrightarrow T \otimes T^*$, and is therefore surjective, as required.

We have now proven that $f_{i} \cS f_{i-1} \cS f_i = f_i \cS f_{i+1} \cS f_i$, for $2 \leq i \leq a-1$.
We therefore have isomorphisms
$$\tilde{A}_i^* \cong {}_i\tilde{T}_{i-1} \underset{\tilde{A}_{i-1}}{\otimes} {}_{i-1}\tilde{T}_i 
\cong {}_i\tilde{T}_{i+1} \underset{\tilde{A}_{i+1}}{\otimes} {}_{i+1}\tilde{T}_i \cong  \tilde{A}_i^*,$$
of $\tilde{A}_i$-$\tilde{A}_i$-bimodules. 
Let us denote this chain of isomorphisms $\phi_i$.
We have 
$$\Hom_{A \otimes A^{op}}(A^*,A^*) \cong \Hom_{A \otimes A^{op}}(A,A) \cong Z(A),$$ and thus $\phi_i$
is multiplication by a central element in $\tilde{A}_i$.
Multiplying the bimodules $_{i}\tilde{T}_{i+1}$ 
by these central elements if necessary, we can assume that in fact $\phi_i = 1$, for $1 \leq i \leq a-1$.

It is now clear that the sum of our bimodule isomorphisms 
$$\tilde{A}_i \cong A_i, \quad {_i\tilde{T}_{i+1}} \cong {_iT_{i+1}}, \quad {_i\tilde{T}^*_{i+1}} \cong {_iT_{i+1}^*}, 
\quad \tilde{A}_i^* \cong A_i^*$$
defines an algebra isomorphism from $\tilde{C}_1^a$ to $C_1^a$, as required.
$\square$

\section{$GL_2$}

In this chapter, we give precise statements of Theorem \ref{Morita} and Conjecture \ref{wishful}, 
together with a justification of theorem \ref{Morita}.

\bigskip 
The determinant representation of $GL_n(F)$ is a polynomial representation of degree $n$.
Therefore, tensoring with the determinant representation defines an exact functor from the 
category of polynomial $GL_n(F)$ representations of degree $r$ to the category of polynomial $GL_n$-representations 
of degree $r+n$, carrying simple modules to simple modules.
Correspondingly, the Schur algebra ${\cS}(n,r)$ can be realised as a quotient 
of ${\cS}(n,r+n)$ by an idempotent ideal ${\cS}(n,r+n)i{\cS}(n,r+n)$.
We denote by ${\cS}(n, \underline{r})$ the inverse limit of the sequence of algebra epimorphisms
$${\cS}(n,r) \twoheadleftarrow {\cS}(n,r+n) \twoheadleftarrow {\cS}(n,r+2n) \twoheadleftarrow...$$
The centre $Z$ of $GL_n(F)$ is isomorphic to $F^\times$, and its group of rational characters is 
therefore isomorphic to $\mathbb{Z}$.  
The category of 
rational representations of $GL_n(F)$ on which $Z$ acts by the character $r \in \mathbb{Z}$ 
is naturally equivalent to $\cS(n, \underline{r}) \ml$.
The category of rational representations of $GL_n(F)$ is therefore isomorphic to the module category of
$\bigoplus_{r \in \mathbb{Z}} \cS(n, \underline{r})$.

For any finite dimensional algebra $A$, the algebra ${\cC}_n(A)$ has an ideal 
$$\bigoplus_{1 \leq i \leq n-1} (A_{i+1} \oplus {_iT_{i+1}} \oplus {_iT_{i+1}^*} \oplus A_i^*),$$ 
the quotient by which is $A_1 \cong A$. In this way, we obtain a sequence of algebra epimorphisms,
$$A \twoheadleftarrow {\cC}_n(A) \twoheadleftarrow {\cC}_n({\cC}_n(A)) \twoheadleftarrow...$$
We denote by $\underleftarrow{\cC}_n(A)$ the inverse limit of this sequence of maps. 
The statement of Conjecture \ref{wishful} is now completely precise:

\bigskip

{\noindent \bf Conjecture \ref{wishful}.} 
\emph{Every block of rational representations of $GL_2(F)$ is equivalent to 
$\underleftarrow{\cC}_p(F) \ml$.}

\bigskip

An equivalent statement is that every block of ${\cS}(2, \underline{r})$ 
is Morita eqivalent to $\underleftarrow{\cC}_p(F)$.
Another is that $\cS \cong \cS_{gr}$, in the notation of the last chapter.

\bigskip

We now give some corollaries of our work in chapter \ref{Schur}. Let $\cS, \cN, A, T$ be as defined there,
and let $U$ be an $\cS$-$\cS$-tilting bimodule.

\begin{lem}
We have $\cN U = U \cN$, and $\cN \cS^* = \cS^* \cN$.
\end{lem}
Proof: 

A tilting bimodule for $\cS$ is given by $U = (\bigoplus_{i=1}^{p-1} \cS f_i) \oplus T$.
Thus, 
$$\cN U = \bigoplus_{1 \leq i,j \leq p-1} f_j \cN f_i = U \cN.$$
We have $\cS^* \cong (\bigoplus_{i=1}^{p-1} \cS f_i) \oplus \cS f_p^*$.
Making this identification, we find  
$$\cN \cS^* \cong (\bigoplus_{1 \leq i,j \leq p-1} f_j \cN f_i) \oplus A_p^* \cong \cS^* \cN. \quad \square$$

\begin{cor}
The space
$$U_{gr} = \bigoplus_{i=0,1,2} \cN^i U/\cN^{i+1} U$$
is a $\cS_{gr}$-$\cS_{gr}$-tilting bimodule.
The space $$(\cS^*)_{gr}  = \bigoplus_{i=0,1,2} \cN^i \cS^*/\cN^{i+1} \cS^*$$
is a $\cS_{gr}$-$\cS_{gr}$-bimodule, isomorphic to $(\cS_{gr})^*$. $\square$
\end{cor}

By Theorem \ref{grMorita}, $\cS_{gr}$ is Morita equivalent to $\cC_p(A) \oplus F^{\oplus I_0}$, 
where $A$ is another Ringel self-dual
Schur algebra, by induction we obtain the following: 

\begin{cor} \label{block}
Then there is a filtration of $\cS$ by ideals, refining the radical filtration, 
whose associated graded ring $\cG$ is Morita equivalent to a direct sum of algebras of the form 
$\underleftarrow{\cC}_p^d(F)$, for $d \in \mathbb{Z}_+$. $\square$
\end{cor}

\bigskip

Given $r \in \mathbb{Z}_+$, we choose $d \geq r$, such that $\cS = \cS(2,d)$ is Ringel self-dual, and $d = r$ (mod $2$).
We have $\cS(2,r) \cong \cS /\cS j\cS$, for some idempotent $j$.
We define $\cG(2,r)$ to be $\cG/\cG j \cG$, where $\cG$ is the graded ring associated to $\cS$ by
Corollary \ref{block}. 
The algebra $\cG(2,r)$ is independent of choice of $d$, and we have algebra epimorphisms
$${\cG}(2,r) \twoheadleftarrow {\cG}(2,r+2) \twoheadleftarrow {\cG}(2,r+4) \twoheadleftarrow...$$
between graded rings ${\cG}(2,r) = {\cG}({\cS}(2,r))$ of Schur algebras.

The statement of Theorem \ref{Morita} is now completely precise. 
Its truth is clear from Corollary \ref{block}. 

\bigskip

{\noindent \bf Theorem \ref{Morita}.} 
\emph{Every block of $\cG(2, \underline{r})$ is Morita equivalent to 
$\underleftarrow{\cC}_p(F)$.} $\square$

\bigskip

\section{Stable equivalence}

A deep conjecture of M. Brou\'e 
predicts that a block of a finite group abelian of abelian defect 
is equivalent to its Brauer correspondent \cite{Broue}.
R. Rouquier has proved that the snag to an inductive proof of this conjecture 
is the lifting of a stable equivalence to a derived equivalence; 
he has also observed an analogy between this and a basic problem in algebraic geometry:
proving that birational Calabi-Yau varieties have equivalent derived categories \cite{Rouquier}, \cite{Bridgeland}.
In this chapter, we prove that the ability to overcome such difficulties would 
also facilitate a proof of Conjecture \ref{wishful}.

We define here a pair of infinite dimensional self-injective quasi-hereditary algebras 
${\cL}_1$ and ${\cL}_2$. 
We define a stable equivalence between ${\cL}_1$ and
${\cL}_2$, sending simple modules to simple modules.
If we could lift this stable equivalence to a Morita equivalence, 
we would have a proof of Conjecture \ref{wishful}.

For background on triangulated categories, we refer to Neeman's book \cite{Neeman}. 
For a concrete approach, 
and a proof that the stable module category of a self-injective algebra is triangulated, 
the reader may consult the book of Happel \cite{Happel}.

\bigskip

Let $\cS, A$ denote the Schur algebras defined in chapter four.
Let $\cQ$ denote the algebra $f \cS f$, where $f= \sum_{i=1}^{p-1} f_i$.

Let ${\cL}_1 = C= C(A)$. 
We now define an algebra ${\cL}_2$, by removing a copy of $C_1^n$ from ${\cL}_1$, 
and gluing a copy of $\cQ$ in its place.

Let $_0T_1$ be a tilting $A$-$\alpha_{1}$-bimodule, and $_{p-1}T_p$ a tilting $\alpha_{p-1}$-$A$-bimodule.
We have canonical bimodule isomorphisms,
$$_0T_1^* \underset{A}{\otimes}{} _0T_1 \cong \alpha_1^*, 
\quad {_0T_1} \underset{\alpha_1}{\otimes} {_0T_1^*} \cong A^*,$$
$${_{p-1}T_{p}^*} \underset{\alpha_{p-1}}{\otimes} {_{p-1}T_{p}} \cong A^*, 
\quad {_{p-1}T_{p}} \underset{A}{\otimes} {_{p-1}T_{p}^*} \cong \alpha_{p-1}^*.$$
We define an algebra $\cL_2$ in the following way:
It consists of three subalgebras $C_{J'_0}$, $\cQ$ and $C_p$ glued together with the bimodules 
$_{0}T_{1}$, $_{0}T_{1}^*$, $_{p-1}T_{p}$, $_{p}T_{p-1}^*$, 
$$\cL_2 := \left(   
 \begin{array}{ccc}
 C_{J'_0} & {}_{0}T_1& 0 \\
 {}_{0}T_{1}^*& \cQ & {}_{p-1}T_{p}\\
 0 & {}_{p-1}T_{p}^*& C_p\\
\end{array}
\right)$$
where the multiplications
$$C_{J'_0} \otimes {}_{0}T_{1} \rightarrow {}_{0}T_{1} \quad {}_{0}T_{1}^* \otimes C_{J'_0} \rightarrow {}_{0}T_{1}^* $$
$${}_{0}T_{1} \otimes \cQ \rightarrow {}_{0}T_{1} \quad \cQ \otimes {}_{0}T_{1}^* \rightarrow {}_{0}T_{1}^*$$
$$C_p \otimes {}_{p-1}T_{p}^*\rightarrow {}_{p-1}T_{p}^* \quad {}_{p-1}T_{p} \otimes C_p \rightarrow {}_{p-1}T_{p}$$
$${}_{p-1}T_{p}^* \otimes  \cQ \rightarrow {}_{p-1}T_{p}^* \quad \cQ \otimes {}_{p-1}T_{p} \rightarrow {}_{p-1}T_{p}$$
are given by the action of the corresponding quotient $A_i$ of the involved algebra $R \in \{C_{J'_0}, \cQ, C_p\}$ 
on the tilting module, the kernel of the surjection $R \twoheadrightarrow A_i$ acting as zero.
The multiplications between the tilting modules
$${}_0T_1 {\otimes}{}_0T_1^* \rightarrow A^*_{0} \subset C_{J'_0},$$
$${}_0T_1^* {\otimes}{}_0T_1 \rightarrow  \alpha^*_{1} \subset {\cQ},$$
$${}_{p-1}T_{p} {\otimes}{_{p-1}T_{p}^*} \rightarrow \alpha^*_{p-1} \subset {\cQ},$$
$${}_{p-1}T_{p}^*{\otimes} {_{p-1}T_{p}} \rightarrow A^*_{p} \subset C_p,$$
are given by the canonical isomorphisms above. 
All other products between elements of the bimodules are zero. 
Similarly multiplying elements of two different subalgebras out of the three yields zero.

\bigskip

Let ${\cL}$ be either ${\cL}_1$ or ${\cL}_2$.
We prove statements concerning both these algebras.

\bigskip

When ${\cL}= {\cL}_2$, we let $f_i$ denote $1_{A_i}$, for $i \in \mathbb{Z}$.

When ${\cL}= {\cL}_2$, we let $A_i$ denote the component $A_i$ of $C_{J'_0}$ for $i \leq 0$, 
the component $\alpha_i$ for $1 \leq i \leq p-1$, and the component $A_i$ of $C_p$ for $i \geq p$.  
We also let $_iT_{i+1}$ denote the component $_iT_{i+1}$ of $C_{J'_0}$ for $i < 0$, 
the component $X_i$ for $1 \leq i \leq p-2$, and the component $T_i$ of $C_p$ for $i \geq p$.
We let $f_i$ denote the idempotent $1_{A_i}$ of $C_{J'_0}$ for $i \leq 0$, 
the idempotent $f_i \in \cS$ for $1 \leq i \leq p-1$, and the idemponent $1_{A_i}$ of $C_p$ for $i \geq p$.

In either case ${\cL} = {\cL}_1, {\cL}_2$, we let $B$ denote that quotient, 
obtained be factoring out the ideal $\oplus_{i \in \mathbb{Z}} ({_iT_{i+1}^*} \oplus A_i^*)$.
We let $B^t$ denote that quotient, 
obtained be factoring out the ideal $\oplus_{i \in \mathbb{Z}} ({_iT_{i+1}} \oplus A_i^*)$.
The algebras $B, B^t$ are both isomorphic to the algebra $B = B(A)$ defined in chapter 3.

Let $\Lambda_{\cL}^1 = \Lambda_B^1$, and $\Lambda_{\cL}^2 = \Lambda_B^2$.

\begin{prop}
The algebra ${\cL}$ is self-injective, 
and quasi-hereditary with respect to both $\Lambda^1_{\cL}$, and 
$\Lambda^2_\cL$.
\end{prop}
\emph{Proof.} 
Note that for any $i$, $\cL_2 f_i$ has the same filtration as $Cf_i$ by a submodule $B^* f_i$ 
and a quotient $Bf_i$. 
Using the same definition of standard modules as for $C$ (which were the same as for $B$) 
we see that by the same arguments as for $C$, 
$\cL_2$ is quasi-hereditary. 
The other quasihereditary structure comes from the filtration of $C1_{A_i}$ 
with quotient $B^tf_i$ and submodule $B^{t*} f_i$.

To show that $\cL_2$ is selfinjective, consider ${}_{\cL_2}\cL_2^*$ which has components
$$\left(   
 \begin{array}{ccc}
 C_{J'_0}^* & {}_{0}T_1 & 0 \\
 {}_{0}T_{1}^* & \cQ^* & {}_{p-1}T_{p} \\
 0 & {}_{p-1}T_{p}^* & C_p^*\\
\end{array} 
\right),$$
where $C_{J'_0}^* = \bigoplus_{i \leq 0} \Hom(C_{J'_0}f_i, F)$, and 
$C_p^* = \bigoplus_{i \geq p} \Hom(C_pf_i, F)$.
Taking into account the self-injectivity of $\cQ$,
which holds since projective modules ${\mathcal Q}f_i, 1 \leq i \leq p-1$ are tilting, 
the selfinjectivity of $C_{J'_0}$ and $C_p$, 
which hold since every projective for $C$ is selfdual and this isn't changed by cutting to a heredity ideal, 
and the fact that $T \cong T^*$ for a tilting $A$-$A$-bimodule $T$, 
we see that this is isomorphic to $\cL_2$ as a left $\cL_2$-module, proving the claim. 
$\square$

\begin{lem} \label{generate}
${\cL}\sml$
is generated by $\add A_p, \add A_{p+1}$.
\end{lem}

\emph{Proof.}
Let ${\mathcal T}$ be the smallest triangulated subcategory of ${\cL}\sml$
containing the subcategories $\add A_p, \add A_{p+1}$.
To prove the lemma, we show by induction on $n$ 
that all $A_i$-modules are in ${\mathcal T}$, for $p+1-n \leq i \leq p+n$.
Since $A$ has finite global dimension, 
we can form a finite projective resolution of length $m$ in $A_i \ml$ for every $A_i$-module ($i=p,p+1$) $M$. 
Then $\Omega^m(M)$ is in $\add A_i$, and by considering triangles 
stemming from exact sequences in $A_i \ml$ we see that $\Omega^{m-1}M, \Omega^{m-2}M, \dots, M \in \mathcal T$. 
Thus the statement holds in case $n=1$.
Suppose the statement is true for $n=N$.
There is a triangle of $\cL$-modules
$$A_i^* \rightarrow \Omega(A_i) \rightarrow {}_{i-1}T_i \oplus {}_{i}T_{i+1}^* \leadsto $$
in ${\mathcal T}$ for $p+1-N \leq i \leq p+N$. Putting $i=p+1-N$, 
we conclude that ${}_{p-N}T_{p+1-N} \in {\mathcal T}$. Putting $i=p+N$, we conclude that
${}_{p+N}T_{p+N+1}^* \in {\mathcal T}$.  
Taking direct summands, we find that $\add {}_{p-N}T_{p+1-N}, \add {}_{p+N}T_{p+N+1}^* \subset {\mathcal T}$. 
Since there is a finite resolution 
$A \hookrightarrow U_1 \rightarrow \cdots \twoheadrightarrow  U_r$ with  $U_j \in \add T$ for all $j$,  
we find that the inductive hypothesis
is true for $n=N+1$, as required. $\square$

\bigskip

Let ${}_A\widehat{T}_A = \dots \rightarrow \widehat T^{(2)} \rightarrow \widehat T^{(1)}\rightarrow \widehat T^{(0)}$ 
be a projective resolution of the bimodule ${}_A T^*_A$ of minimal length. 
This is finite since $A \otimes A^{op}$ is finite-dimensional 
and quasi-hereditary, and therefore of finite global dimension.
Then the total complex of ${}_A\widehat{T}_A  \underset{A}{\otimes} \cdots \underset{A}{\otimes}{}_A\widehat{T}_A$ 
($k$ factors) is a projective resolution of 
${}_A{T}^*_A  \underset{A}{\otimes} \cdots \underset{A}{\otimes}{}_A{T}^*_A$ ($k$ factors).
Also, for any primitive idempotent $e$ of $A$,  ${}_A\widehat{T}e$ is a projective resolution 
of the indecomposable summand $T^*e$ of the tilting module.

Let $q: \cL \twoheadrightarrow B$ be the algebra epimorphism of ${\cL}$ onto $B$, which factors out the ideal
$\bigoplus_{i \in \mathbb{Z}}( _iT_{i+1}^* \oplus A_i^*)$.
The endomorphism ring of $Bf_i$ is $A_i$.
Let ${\cK} = ker(q)$.
We have an isomorphism ${\cK} f_i \cong Bf_{i+1} \underset{A_{i+1}}{\otimes}{}_{i}T_{i+1}^*$ of $(\cL, A_i)$-bimodules,
coming from the right multiplication with ${}_{i}T_{i+1}$. 

Given a complex of $A_i$-modules $C$, 
we define ${\cL}f_i \underset{A_i}{\odot} C$ to be some complex of projective
${\cL}$-modules such that $A_i \underset{\cL}{\otimes} {\cL}f_i \underset{A_i}{\odot} C \cong C$.
Such a complex exists, by the usual lifting argument, but is not necessarily unique. 
However, ${\cL}f_i \underset{A_i}{\odot} C$ contains a canonical subcomplex 
isomorphic to ${\cK}f_i \underset{A_i}{\otimes} C$, 
the quotient by which is isomorphic to ${Bf_i} \underset{A_i}{\otimes} C$.

\begin{lem}\label{projdiminfty}
Every direct summand of ${}_{\cL}Bf_i$ has infinite projective dimension.
\end{lem}

\emph{Proof}
Let $e$ be any idempotent which is a summand of $f_i$.
We manufacture an infinite projective resolution of $_{\cL} Be$, 
whose components in degree $>>0$ are components of $\cL f_j$, for $j>>0$. 
This allows us to prove that given any $m>>0$, 
there exists a simple $A_j$-module $L$, for some $j>>i$, 
such that $Ext^m(Be, L) \neq 0$.

We have an isomorphism ${\cK} f_i \cong Bf_{i+1} \underset{A_{i+1}}{\otimes}{}_{i}T_{i+1}^*$.
Right multiplication by $e$ gives an isomorphism of $\cL$-modules 
${\cK} e \cong Bf_{i+1} \underset{A_{i+1}}{\otimes}{}_{i}T_{i+1}^*e$.  
So we have short exact sequences
\begin{equation}\label{resstart}Bf_{i+1} \underset{A_{i+1}}{\otimes}{}_{i}T_{i+1}^* \rightarrow 
\cL f_i \rightarrow Bf_i\end{equation}
and
\begin{equation}\label{resstarte}Bf_{i+1} \underset{A_{i+1}}{\otimes}{}_{i}T_{i+1}^*e \rightarrow 
\cL e \rightarrow Be.\end{equation}
We thus obtain a natural map of complexes of $\cL$-modules,
$$\phi_i: Bf_{i+1} \underset{A_{i+1}}{\otimes} {_{i+1}\widehat{T}_i}  \rightarrow {\cL}f_i$$
with direct summand
$$\phi_e: Bf_{i+1} \underset{A_{i+1}}{\otimes} {_{i+1}\widehat{T}_i}e  \rightarrow {\cL}e.$$
Since all modules occurring in the complexes ${}_{i+1}\widehat{T}_{i}$ are projective as left $A_{i+1}$-modules,
and hence $\Delta$-filtered, and 
$$Bf_{i+1} \underset{A_{i+1}}{\otimes} -  \cong (A_{i+1} \oplus {_iT_{i+1}}) \underset{A_{i+1}}{\otimes} -$$ 
is exact on $\mathcal{F}(\Delta)$, 
we find that $Bf_{i+1} \underset{A_{i+1}}{\otimes} {}_{i+1}\widehat{T}_{i}$ is quasi-isomorphic to 
$Bf_{i+1} \underset{A_{i+1}}{\otimes}{}_{i+1}T_{i}$, 
and $Bf_{i+1} \underset{A_{i+1}}{\otimes} {}_{i+1}\widehat{T}_ie$ is quasi-isomorphic to 
$Bf_{i+1} \underset{A_{i+1}}{\otimes}{}_{i}T_{i+1}e$ . 
Therefore, the cone of $\phi_i$ is quasi-isomorphic to $Bf_i$ and the cone of $\phi_e$ is quasi-isomorphic to $Be$.

We have an exact sequence of complexes of $\cL$-modules
$$Bf_{i+2} \underset{A_{i+2}}{\otimes}{}_{i+2}T_{i+1} \underset{A_{i+1}}{\otimes}{}_{i+1} \widehat T_{i}\rightarrow 
\cL f_{i+1} \underset{A_{i+1}}{\odot}{}_{i+1}\widehat T_{i}\rightarrow 
Bf_{i+1} \underset{A_{i+1}}{\otimes}{}_{i+1}\widehat T_{i},$$
again with direct summand 
$$Bf_{i+2} \underset{A_{i+2}}{\otimes}{}_{i+2}T_{i+1} \underset{A_{i+1}}{\otimes}{}_{i+1} \widehat T_ie\rightarrow 
\cL f_{i+1} \underset{A_{i+1}}{\odot}{}_{i+1}\widehat T_ie\rightarrow 
Bf_{i+1} \underset{A_{i+1}}{\otimes}{}_{i+1}\widehat T_ie.$$
These are exact in every degree, since the sequence in a given degree is
obtained by tensoring the short exact sequence (\ref{resstart}) 
with the corresponding entry of ${}_{i+1}\widehat T_{i}$ (resp.\ ${}_{i+1}\widehat T_ie$), 
which is flat as a left $A_{i+1}$-module.
Therefore $Bf_{i+1} \underset{A_{i+1}}{\otimes}{}_{i+1}\widehat T_ie$
is quasi-isomorphic to the total complex of 
$$Bf_{i+2} \underset{A_{i+2}}{\otimes}{}_{i+2}T_{i+1}
\underset{A_{i+1}}{\otimes}{}_{i+1} \widehat T_{i}\rightarrow \cL
f_{i+1} \underset{A_{i+1}}{\odot}{}_{i+1}\widehat T_ie.$$

We now claim that $Bf_{i+2} \underset{A_{i+2}}{\otimes}{}_{i+2}T_{i+1}
\underset{A_{i+1}}{\otimes}{}_{i+1} \widehat T_ie$ is quasi-isomorphic
to the total complex of $Bf_{i+2}
\underset{A_{i+2}}{\otimes}{}_{i+2}\widehat T_{i+1}
\underset{A_{i+1}}{\otimes}{}_{i+1} \widehat T_ie$. Indeed, $Bf_{i+2}
\underset{A_{i+2}}{\otimes}{}_{i+2}T_{i+1}$ is quasi-isomorphic to
$Bf_{i+2} \underset{A_{i+2}}{\otimes}{}_{i+2}\widehat T_{i+1}$ by the
above, and since every module occurring in ${}_{i+1} \widehat T_ie$ is
projective as a left $A_{i+1}$-module and every module occurring in
$Bf_{i+2} \underset{A_{i+2}}{\otimes}{}_{i+2}\widehat T_{i+1}$ is
projective as a right $A_{i+1}$-module, the rows and columns in the
double complex
%
%
are exact, proving the claim.

We now know that $Be$ is quasi-isomorphic to the total complex of
$$Bf_{i+2}
\underset{A_{i+2}}{\otimes}{}_{i+2}\widehat T_{i+1} \underset{A_{i+1}}{\otimes}{}_{i+1} \widehat T_ie \rightarrow 
{\cL}f_{i+1} \underset{A_{i+1}}{\odot}{}_{i+1} \widehat T_ie \rightarrow {\cL}e.$$
Iterating this procedure, we obtain
a projective resolution $P_.(Be)$ of ${}_{\cL} Be$,
with a filtration whose sections are isomorphic to 
$$\cL f_j \underset{A_j}{\odot}{}_j \widehat{T}_{j-1} \underset{A_{j-1}}{\otimes} {}_{j-1}\widehat{T}_{j-2} 
\underset{A_{j-2}}{\otimes} \cdots \underset{A_{i+1}}{\otimes}{}_{i+1} \widehat{T}_ie [i-j],$$
as $\cL$-modules for $j \geq i$.

We now claim that for every $m \geq 0$, there exists an irreducible $\cL$-module $L$ such that 
$\Ext^m_\cL(Be, L)\neq 0$.
From the projective resolution above we see that 
$$P_m(Be) \cong$$ 
$$\left( \cL f_{i+m} \odot (\widehat T^{(0)})^{\otimes m}e \right) \oplus 
\left( \underset{\substack{k<m \\\sum_{j=1}^k r_j =m-k }}{\bigoplus} 
\cL f_{i+k} \odot \widehat T^{(r_1)} \otimes \cdots \otimes \widehat T^{(r_k)}e \right),$$
where all tensors are taken over $A_i$, for some $i$.
Now choose an irreducible $A_{i+m}$-module $L$ in the head of $(\cL f_{i+m} \odot (\widehat T^{(0)})^{\otimes m})e$.
Then $\Hom_\cL (P_m(Be), L) \neq 0$ but $\Hom_\cL (P_k(Be), L) = 0$ for all $k<m$.
Furthermore $$\Hom_\cL (P_{m+1}(Be), L) =$$ 
$$\overset{m-1}{\underset{r=0}{\bigoplus}}
\Hom_\cL (\cL f_{i+m} \odot (\widehat T^{(0)})^{\otimes r} \otimes \widehat T^{(1)} 
\otimes (\widehat T^{(0)})^{\otimes m-r-1}e, L).$$
The map $$\overset{m-1}{\underset{r=0}{\bigoplus}}\cL f_{i+m} \odot (\widehat T^{(0)})^{\otimes r} \otimes 
\widehat T^{(1)} \otimes (\widehat T^{(0)})^{\otimes m-r-1}e \rightarrow 
\cL f_{i+m} \otimes (\widehat T^{(0)})^{\otimes m}e,$$
is not surjective as the cokernel at least has a quotient $\cL f_{i+m} \odot (T^{\otimes m})e \neq 0$.
Therefore we can choose $L$ such that not every morphism in $\Hom_\cL (P_m(Be), L)$ 
comes from a morphism in $\Hom_\cL (P_{m+1}(Be), L)$. Hence $\Ext^m_\cL(Be, L)\neq 0$.
$\square$

\bigskip

Let $\cL^+ = \cL^{I_p'}$. Thus, $\cL^+$ is a quasi-hereditary quotient of $\cL$ with poset $I_p'$. 
By definition, there is an isomorphism $\cL_1^+ \cong \cL_2^+$.

By the theory of E. Cline, B. Parshall, and L. Scott (\cite{CPS:88}, Theorem 3.9), 
we have an embedding of derived categories
$$j: D^b(\cL^+) \rightarrow D^b(\cL).$$

By a theorem of Rickard \cite{R:89}, we have a Verdier quotient of triangulated categories,
$$\pi: D^b(\cL) \rightarrow \cL \sml,$$
whose kernel is the thick subcategory of perfect complexes. Note that
even though his theorem only includes finite-dimensional
self-injective algebras, the same proof goes through in the locally
finite-dimensional case.

\begin{prop} \label{quotient}
The composition
$$D^b(\cL^+) \overset{j}{\rightarrow} D^b(\cL) \overset{\pi}{\rightarrow} \cL\sml$$
has dense image,
and kernel ${\mathcal T}$, where ${\mathcal T}$ is the thick subcategory of $D^b(\cL^+)$
generated by $\{ \cL^+f_i, i>p \}$.
\end{prop}
\emph{Proof.} In the above composition of functors, $M \in \cL^+ \ml \subset D^b(\cL^+) $ 
maps to the isomorphism class of $M$ in $\cL \sml$. 
Considering $\add A_p, \add A_{p+1} \subset  \cL^+ \ml$, and applying Lemma \ref{generate}, 
we see that the image is indeed dense.

It is obvious that $\mathcal T$ is contained in the kernel, since for $i<0$, 
the projectives for $\cL 1_{A_i}$ and $\cL^+ 1_{A_i}$ are the same, 
so under the inclusion $j$ bounded complexes in projectives from $\cL^+ 1_{A_i}$ ($i<0$) map 
to bounded complexes in projectives for $\cL$, which become isomorphic to zero under $\pi$.

Suppose that $E_.$ is a bounded complex of projective modules in $K^b({\cL}^+)$ of minimal length, 
such that $E_. \notin {\mathcal T}$.
Therefore, some direct summand of $\cL 1_{A_0}$ occurs in $E_.$. 
By cutting in the ``stupid way'' and shifting in degree
we may assume that $E_0 \neq 0$ is a direct sum of summands of $\cL^+ 1_{A_0}$,
and $E_i = 0$, for $i>0$. 
The image of $E_.$ under $j$ in $D^b(\cL) \cong K^{-,b}(\cL \projl)$ does not have a presentation 
as a finite complex of projective $\cL$-modules by Lemma \ref{projdiminfty}, 
and therefore $E_.$ is not contained in the kernel of $\pi$.
This completes the proof of the proposition. $\square$

\begin{thm} \label{stable}
We have a stable equivalence,
$$\cL_1 \sml \cong \cL_2 \sml,$$
sending simple modules to simple modules.
\end{thm}
Proof:

We have $\cL_1^+ \cong \cL_2^+$. Therefore, the
stable equivalence is immediate from proposition \ref{quotient}.
The fact that simple modules correspond to simple modules is obvious on the subcategory of $\cL_1^+$-modules.
For the remaining simples, one proceeds by induction, exactly as in the proof of Lemma \ref{generate}.
$\square$

\section{Epilogue}

We end with some remarks and open questions.

\begin{remark}
If there is a stable equivalence of Morita type between two finite dimensional algebras, sending simple modules to 
simple modules, then a theorem of M. Linckelmann states that the algebras are in fact Morita equivalent \cite{L:96}.
If there is a stable equivalence of Morita type between two finite dimensional algebras, one of which is graded, then
a theorem of R. Rouquier states that there is a compatible grading on the other algebra \cite{Rouquier00}.
However, we are unable to apply these results, 
since the stable equivalence of theorem~\ref{stable} is not manifestly of Morita type. 

Is it the case that any stable equivalence between locally finite dimensional algebras, 
one of which has a grading refining the radical filtration, and which sends simple modules to simple modules, 
must lift to a Morita equivalence ? If this were so, then Conjecture \ref{wishful} would follow.
\end{remark} 

\begin{remark}
The problem of finding gradings on modular representation categories is rather a general one,
related to the celebrated conjecture of G. Lusztig concerning irreducible characters of algebraic groups 
(see \cite{Jantzen}).
For example, one expects blocks of Schur algebras ${\cS}(n,n)$ to have a grading refining the radical filtration, 
at least when the weight of the block is less than $p$.

We have conjectured that blocks of Schur algebras ${\cS}(n,n)$ 
are all derived equivalent to certain subquotients 
of a symmetric quasi-hereditary algebra, the Schiver double ${\mathcal D}_{A_\infty}$ (see \cite{Turner}, \cite{Turner2}). 
The most obvious barrier to a proof of this is the difficulty of finding a 
grading on the Rock blocks.
Conjecture \ref{wishful} can be thought of as a simple analogue of the Schiver double conjecture,
the algebra ${\cC}_p(A)$ playing a similar role in this paper, 
to that played by the algebra ${\mathcal D}_{A_\infty}$ in the theory of Rock blocks.
Indeed, the development of Conjecture \ref{wishful} was made, 
with a view towards understanding the Schiver double conjecture better.
We hope the method of defining stable equivalences introduced in this paper may prove useful, 
as a step towards a proof of the Schiver double conjecture.  
\end{remark}

\begin{remark}
It would be interesting if there were analogues of Conjecture \ref{wishful} 
for algebraic groups other than $GL_2$. 
Let us speculate on what features such generalisations might possess.

Suppose that $T$ is a tilting bimodule for $A$.
Then we may think of the pair $(A,T)$ as defining a collection  
of triangulated categories and exact functors,
$$... \xymatrix{
D^b(A \ml) \ar@/{}^{-.7pc}/[r]_{\alpha_1} & D^b(A \ml) \ar@/{}_{.7pc}/[l]_{\beta_1}
\ar@/{}^{-.7pc}/[r]_{\alpha_2} & D^b(A \ml) \ar@/{}_{.7pc}/[l]_{\beta_2}
\ar@/{}^{-.7pc}/[r]_{\alpha_3} & D^b(A \ml) \ar@/{}_{.7pc}/[l]_{\beta_3} } ...$$
such that $\alpha_i \beta_i \cong \beta_{i-1} \alpha_{i-1}$, for $i \in \mathbb{Z}$.
Here, $\alpha_i = \beta_i = - \otimes_A^L T$.
Passing to Grothendieck groups, we obtain a free representation of
the preprojective algebra $\Pi_\infty$ on a quiver orienting an infinite line.
In other words, the pair $(A,T)$ defines a $\Pi_\infty$-category, which we denote ${\mathcal F}(A,T)$ 
(cf. \cite{BFK}, \cite{ChRou}).
The map ${\mathcal C}_n$ can therefore be thought of as a map 
$${\cC}_n \circlearrowright \Pi_\infty \cat,$$
where $\Pi_\infty \cat$ denotes a collection of $\Pi_\infty$-categories, 
taking ${\mathcal F}(A,T)$ to ${\mathcal F}({\mathcal C}_n(A,T))$. 
Passing to Grothendieck groups, we see that 
$$K({\mathcal F}({\mathcal C}_n(A,T))) \cong K({\mathcal F}(A,T))^{\oplus n}.$$ 
One way to think of the map ${\mathcal C}_n$ is therefore as a categorification of the 
functor $-^{\oplus n}$ on $\Pi_\infty \ml$. 
To be more precise, one should define $\Pi_\infty \cat$ as a $2$-category,
and $\cC_n$ as an endo-$2$-functor of $\Pi_\infty \cat$. 

Let $q \in \mathbb{C}^\times$ be a $p^{th}$ root of unity. 
Any block of the quantum group $q$-$GL_{2}(\mathbb{C})$ is Morita equivalent to
$C_1(\mathbb{C})$, the Brauer tree algebra on a semi-infinite line 
(\cite{KX98}, Corollary 7.3).
Translation by $p$ embeds a semi-infinite line in itself, and we have a corresponding 
algebra monomorphism from $C_1(\mathbb{C})$ to itelf, 
related to Steinberg's tensor product theorem on $q$-$GL_2(\mathbb{C})$. 
By composition, we obtain a sequence of embeddings,
$$C_1(\mathbb{C}) \hookrightarrow C_1(\mathbb{C}) \hookrightarrow C_1(\mathbb{C}) \hookrightarrow ...,$$
whose direct limit is $C(\mathbb{C})$, the Brauer tree algebra on an infinite line.
The preprojective algebra $\Pi_\infty$ is the Koszul dual of $C(\mathbb{C})$.

To explore the modular representation theory of blocks of a reductive algebraic group $G(F)$,
as we have done in this paper in case $G=GL_2$, 
one should perhaps first look for Koszul duals 
of direct limits of blocks of the corresponding quantum group 
$q$-$G(\mathbb{C})$, 
before looking to define categories over these Koszul duals,
and categorifications of functors between their module categories.
\end{remark}

\bibliographystyle{amsplain}

\begin{thebibliography}{10}
\bibitem{BFK}
Joseph Bernstein, Igor Frenkel, and Mikhail Khovanov, \emph{A categorification
  of the {T}emperley-{L}ieb algebra and {S}chur quotients of
  {$U(\mathfrak{sl}\sb 2)$} via projective and {Z}uckerman functors}, Selecta
  Math. (N.S.) \textbf{5} (1999), no.~2, 199--241. \MR{MR1714141 (2000i:17009)}

\bibitem{Bridgeland}
Tom Bridgeland, \emph{Derived categories of coherent sheaves},  (2006), ICM
  survey article, available at {http://arxiv.org/abs/math.AG/0602129}.

\bibitem{Broue}
Michel Brou{\'e}, \emph{Isom\'etries parfaites, types de blocs, cat\'egories
  d\'eriv\'ees}, Ast\'erisque (1990), no.~181-182, 61--92. \MR{MR1051243
  (91i:20006)}

\bibitem{ChRou}
Joseph Chuang and Rapha\"el Rouquier, \emph{Derived equivalences for symmetric
  groups and $\mathfrak {sl}_2$-categorification},  (2004),
  {http://www.maths.leeds.ac.uk/~rouquier/papers.html}.

\bibitem{CPS:88}
Edward Cline, Brian Parshall, and Leonard Scott, \emph{Finite-dimensional
  algebras and highest weight categories}, J. Reine Angew. Math. \textbf{391}
  (1988), 85--99. \MR{MR961165 (90d:18005)}

\bibitem{Dlab}
Vlastimil Dlab, \emph{Quasi-hereditary algebras revisited}, An. \c Stiin\c t.
  Univ. Ovidius Constan\c ta Ser. Mat. \textbf{4} (1996), no.~2, 43--54,
  Representation theory of groups, algebras, and orders (Constan\c ta, 1995).

\bibitem{D:98}
Stephen Donkin, \emph{The {$q$}-{S}chur algebra}, London Mathematical Society
  Lecture Note Series, vol. 253, Cambridge University Press, Cambridge, 1998.
  \MR{MR1707336 (2001h:20072)}

\bibitem{EH:02}
Karin Erdmann and Anne Henke, \emph{On {R}ingel duality for {S}chur algebras},
  Math. Proc. Cambridge Philos. Soc. \textbf{132} (2002), no.~1, 97--116.
  \MR{MR1866327 (2002j:20081)}

\bibitem{Green}
James~A. Green, \emph{Polynomial representations of {${\rm GL}\sb{n}$}},
  Algebra, Carbondale 1980 (Proc. Conf., Southern Illinois Univ., Carbondale,
  Ill., 1980), Lecture Notes in Math., vol. 848, Springer, Berlin, 1981,
  pp.~124--140. \MR{MR613180 (82i:20019)}

\bibitem{Happel}
Dieter Happel, \emph{Triangulated categories in the representation theory of
  finite dimensional algebras}, London Math Soc. Lecture Note Ser., vol. 119,
  Cambridge Univ. Press, Cambridge, 1988.

\bibitem{HK:02}
Anne Henke and Steffen Koenig, \emph{Relating polynomial {${\rm
  GL}(n)$}-representations in different degrees}, J. Reine Angew. Math.
  \textbf{551} (2002), 219--235. \MR{MR1932179 (2003k:20066)}

\bibitem{Jantzen}
Jens~Carsten Jantzen, \emph{Representations of algebraic groups}, second ed.,
  Mathematical Surveys and Monographs, vol. 107, American Mathematical Society,
  Providence, RI, 2003. \MR{MR2015057 (2004h:20061)}

\bibitem{KX98}
Steffen K{\"o}nig and Changchang Xi, \emph{Strong symmetry defined by twisting
  modules, applied to quasi-hereditary algebras with triangular decomposition
  and vanishing radical cube}, Comm. Math. Phys. \textbf{197} (1998), no.~2,
  427--441. \MR{MR1652759 (99h:16028)}

\bibitem{L:96}
Markus Linckelmann, \emph{Stable equivalences of {M}orita type for
  self-injective algebras and {$p$}-groups}, Math. Z. \textbf{223} (1996),
  no.~1, 87--100. \MR{MR1408864 (97j:20011)}

\bibitem{Neeman}
Amnon Neeman, \emph{Triangulated categories}, Annals of Mathematics Studies,
  vol. 148, Princeton University Press, Princeton, NJ, 2001. \MR{MR1812507
  (2001k:18010)}

\bibitem{ParshallScott}
Brian Parshall and Leonard Scott, \emph{Derived categories, quasi-hereditary
  algebras, and algebraic groups},  (1988),
  http://www.math.virginia.edu/~lls2l/reprnt.htm.

\bibitem{R:89}
Jeremy Rickard, \emph{Derived categories and stable equivalence}, J. Pure Appl.
  Algebra \textbf{61} (1989), no.~3, 303--317. \MR{MR1027750 (91a:16004)}

\bibitem{Ringel}
Claus~Michael Ringel, \emph{The category of modules with good filtrations over
  a quasi-hereditary algebra has almost split sequences}, Math. Z. \textbf{208}
  (1991), no.~2, 209--223. \MR{MR1128706 (93c:16010)}

\bibitem{Rouquier00}
Rapha\"el Rouquier, \emph{Groupes d'automorphismes et \'equivalences stables ou
  deriv\'ees},  (2000), {http://www.maths.leeds.ac.uk/~rouquier/papers.html}.

\bibitem{Rouquier}
\bysame, \emph{Derived equivalences and finite dimensional algebras},  (2006),
  ICM survey article, available at
  {http://www.maths.leeds.ac.uk/~rouquier/papers.html}.

\bibitem{Turner}
Will Turner, \emph{Rock blocks},  (2004),
  {http://www.maths.ox.ac.uk/~turnerw/}.

\bibitem{Turner2}
\bysame, \emph{Tilting equivalences: from hereditary algebras to symmetric
  groups},  (2006), {http://www.maths.ox.ac.uk/~turnerw/}.

\end{thebibliography}

\bigskip

\noindent Vanessa Miemietz,
Universit\"at zu K\"oln,
Mathematisches Institut,
Weyertal 86-90,
D-50931 K\"oln,
Germany.

\bigskip

\noindent Will Turner,
Department of Mathematics, University of Oxford,
Oxford, England.
Email: turnerw@maths.ox.ac.uk

\end{document}